\documentclass[11pt]{amsart}

\usepackage{amssymb, geometry} 
\usepackage[all]{xy}

\newif\ifshow
\showtrue

\newcommand{\ed}{

\section{Concluding remarks}

Most of the results provided here for complete groups, have natural extensions to
incomplete groups. For these extensions, one needs to consider the dual group $\Gh$ with
$[P,\epsilon]$ a neighborhood of the identity for each \emph{precompact} $P\sub G$.
The extension is sometimes straightforward, using Theorem \ref{catchcompact}.

Similarly, the results of Section \ref{AXsec} extend to completely regular spaces that
are not $\mu$-spaces. Here, one should consider \emph{functionally bounded} subsets of
$X$ instead of compact subsets of $X$, and the topology of $C(X,\bbT)$ should be the
functionally bounded-open topology. The main result of this section would then deal with
spaces $X$ having a cofinal family of functionally bounded sets, and whose topology is
determined by its functionally bounded sets. We point out that in this case, the $\mu$-completion
of $X$ is $k_\w$, and $X$ is dense in this completion.

With some adaptation, the results presented here for $k_\w$ groups also apply
to locally convex vector spaces that have a countable cofinal family of bounded sets. For instance, any
countable inductive limit of DF-spaces.

The present work is not the only one where pcf theory arises naturally in a study of a seemingly
unrelated concept. Another recent example is in Feng and Gartside's paper \cite{FengGartside},
where pcf theory turned out essential in a study of a problem motivated by Hilbert's 13th problem.

\subsection*{Acknowledgments}
The second and third named authors 
acknowledge partial support by Generalitat Valenciana,
grant PROMETEO/2014/062, and by Universitat Jaume I, grant P11B2012-05.
We thank Maria Jesus Chasco, Jorge Galindo, and Mikhail Tkachenko
for useful discussions on topological groups,
Moti Gitik and Assaf Rinot for useful discussions on pcf theory,
and Lydia Au\ss{}enhofer and Adi Jarden for reading parts of this paper and making comments.
Some of the work on this paper was carried out when the fourth named
author was visiting the other authors at Universitat Jaume I, Castell\'{o}n.
This author thanks his hosts for their kind hospitality.
We are especially grateful to the referee for a thorough report that helped considerably in
improving the presentation of this paper.

\end{document}
}

\newcommand{\IndexPrint}[1]{#1}
\newcommand{\Ghh}{{G\hat{~}\hat{~}}}
\newcommand{\Gh}{\widehat{G}}
\newcommand{\weight}{\op{w}}
\newcommand{\density}{\op{d}}
\newcommand{\localdensity}{\op{ld}}
\newcommand{\boundedness}{\op{b}}
\newcommand{\Nh}{\widehat{\cN}}
\newcommand{\la}{\langle}
\newcommand{\ra}{\rangle}

\long\def\forget#1\forgotten{}
\newcommand{\cl}[1]{\overline{#1}}
\newcommand{\inv}{^{-1}}
\newcommand{\kw}{\op{kw}}
\newcommand{\w}{\omega}
\newcommand{\Fin}[1]{\op{Fin}(#1)} 
\newcommand{\finw}[1]{{\Fin{#1}^\N}}
\newcommand{\kfinw}{\finw{\kappa}}
\newcommand{\fb}{\mathfrak{b}}
\newcommand{\fd}{\mathfrak{d}}
\newcommand{\fc}{\mathfrak{c}}

\newcommand{\mult}{*}

\newcommand{\alephes}{{\aleph_0}}

\newcommand{\beq}{\begin{eqnarray*}}
\newcommand{\eeq}{\end{eqnarray*}}
\newcommand{\op}[1]{\operatorname{#1}}
\newcommand{\PK}{\op{PK}}
\newcommand{\K}{\op{K}}

\newcommand{\cof}{\op{cof}}

\newcommand{\vphi}{\varphi}
\newcommand{\Bdd}{\op{Bdd}}
\newcommand{\card}[1]{\left|#1\right|}

\newcommand{\Impl}{\Rightarrow}
\newcommand{\bi}{\begin{itemize}}
\newcommand{\itm}{\item}
\newcommand{\ei}{\end{itemize}}
\newcommand{\be}{\begin{enumerate}}
\newcommand{\ee}{\end{enumerate}}

\newcommand{\Union}{\bigcup}
\newcommand{\sub}{\subseteq}
\newcommand{\spst}{\supseteq}
\newcommand{\N}{\mathbb{N}}
\newcommand{\bbQ}{\mathbb{Q}}
\newcommand{\NN}{{\N^\N}}
\newcommand{\Z}{\mathbb{Z}}
\newcommand{\bbC}{\mathbb{C}}
\newcommand{\bbT}{\mathbb{T}}

\newcommand{\R}{\mathbb{R}}
\newcommand{\set}[2]{\{\, #1 : #2\,\}}

\newcommand{\sm}{\setminus}

\newcommand{\cK}{\mathcal{K}}
\newcommand{\cF}{\mathcal{F}}

\newcommand{\cS}{\mathcal{S}}

\newcommand{\cN}{\mathcal{N}}
\newcommand{\x}{\times}

\newtheorem{thm}{Theorem}[section]
\newcommand{\bthm}{\begin{thm}} \newcommand{\ethm}{\end{thm}}
\newtheorem{prop}[thm]{Proposition}
\newcommand{\bprp}{\begin{prop}} \newcommand{\eprp}{\end{prop}}
\newtheorem{fact}[thm]{Fact}
\newcommand{\bfct}{\begin{fact}} \newcommand{\efct}{\end{fact}}
\newtheorem{prob}[thm]{Problem}
\newcommand{\bprb}{\begin{prob}} \newcommand{\eprb}{\end{prob}}
\newtheorem{lem}[thm]{Lemma}
\newcommand{\blem}{\begin{lem}} \newcommand{\elem}{\end{lem}}
\newtheorem{claim}[thm]{Claim}
\newcommand{\bclm}{\begin{claim}} \newcommand{\eclm}{\end{claim}}
\newtheorem{cor}[thm]{Corollary}
\newcommand{\bcor}{\begin{cor}} \newcommand{\ecor}{\end{cor}}
\newtheorem{conj}[thm]{Conjecture}
\newcommand{\bcnj}{\begin{conj}} \newcommand{\ecnj}{\end{conj}}
\theoremstyle{definition}
\newtheorem{defn}[thm]{Definition}
\newcommand{\bdfn}{\begin{defn}} \newcommand{\edfn}{\end{defn}}
\newtheorem{spec}[thm]{Specializing}
\newcommand{\bspc}{\begin{spec}} \newcommand{\espc}{\end{spec}}
\theoremstyle{remark}
\newtheorem{rem}[thm]{Remark}
\newcommand{\brem}{\begin{rem}} \newcommand{\erem}{\end{rem}}
\newtheorem{cnv}[thm]{Convention}
\newcommand{\bcnv}{\begin{cnv}} \newcommand{\ecnv}{\end{cnv}}
\newtheorem{exam}[thm]{Example}
\newcommand{\bexm}{\begin{exam}} \newcommand{\eexm}{\end{exam}}
\newcommand{\bpf}{\begin{proof}} \newcommand{\epf}{\end{proof}}

\title[The character of topological groups]
{The character of topological groups, via
bounded systems, Pontryagin--van Kampen duality and pcf theory}

\author[C. Chis]{Cristina Chis}
\address[Chis]{Universitat Jaume I, Departamento de Mate\-m\'{a}ticas,
Campus de Riu Sec, 12071 Castell\'{o}n, Spain.}
\email{chis@mat.uji.es}
\author[M. V. Ferrer]{M. Vincenta Ferrer}
\address[Ferrer]{Universitat Jaume I, IMAC and Departamento de Mate\-m\'{a}ticas,
Campus de Riu Sec, 12071 Castell\'{o}n, Spain.}
\email{mferrer@mat.uji.es}
\author[S. Hern\'andez]{Salvador Hern\'andez}
\address[Hern\'andez]{Universitat Jaume I, INIT and Departamento de Mate\-m\'{a}ticas,
Campus de Riu Sec, 12071 Castell\'{o}n, Spain.}
\email{hernande@mat.uji.es}
\author[B. Tsaban]{Boaz Tsaban}
\address[Tsaban]{Department of Mathematics, Bar-Ilan University, Ramat Gan 52900, Israel}
\email{tsaban@math.biu.ac.il}



\makeindex

\begin{document}

\begin{abstract}
The Birkhoff--Kakutani Theorem asserts that a topological
group is metrizable if and only if it has countable character.
We develop and apply tools for the estimation of the character
for a wide class of nonmetrizable topological groups.

We consider abelian groups whose topology is determined by a countable cofinal family of
compact sets.
These are the closed subgroups of Pontryagin--van Kampen duals of \emph{metrizable}
abelian groups,
or equivalently, complete abelian groups whose dual is metrizable.
By investigating these connections, we show that also in these cases,
the character can be estimated, and that it is determined by the weights of
the \emph{compact} subsets of the group, or of quotients of the group by compact subgroups.
It follows, for example, that the density and the local density of an abelian
metrizable group determine
the character of its dual group.
Our main result applies to the more general case of closed subgroups of
Pontryagin--van Kampen duals of abelian \v{C}ech-complete groups.

In the special case of free abelian topological groups, our results extend
a number of results of Nickolas and Tkachenko, which were proved using
combinatorial methods.

In order to obtain concrete estimations, we establish a natural bridge between the studied concepts
and pcf theory, that allows the direct application of several major results from that theory.
We include an introduction to these results and their use.
\end{abstract}

\maketitle


\section{Overview and main results}\label{sec:overview}

The topological structure of a topological group is completely determined by its
local structure at an element.
The most fundamental invariant of the local structure is the
\emph{character}\index{character!the character of a topological group}, the minimal cardinality of a local basis.
Metrizable groups have countable character, and the celebrated Birkhoff--Kakutani Theorem asserts
that
this is the only case where the character is countable.

The computation of the character of nonmetrizable groups may be a difficult task. For example,
the character of free abelian topological groups is only known in
some cases (cf. \cite{NT1,NT2}).
The \emph{free abelian topological group $A(X)$} over a Tychonoff space $X$
is the abelian topological group with the universal property that each continuous function $\vphi$ from
$X$ into any abelian topological group $H$ has a unique extension to a continuous homomorphism
$\tilde\vphi\colon A(X)\to H$.
$$\xymatrix{
A(X)\ar^{\exists ! \tilde \varphi}[rd]\\
X~ \ar@^{(->}^{\mathrm{id}}[u]\ar^{\forall \varphi\ }[r] & H
}$$
As a set, $A(X)$ is the family of all formal linear combinations of elements of $X$ over the integers.
But the topology of $A(X)$ is very complex, and in general, it is not known
how to determine the character of $A(X)$ from the properties of $X$.

In this paper, we make use of the fact that groups from an important class of topological groups,
whose character estimation was intractable for earlier methods,
contain open subgroups whose Pontryagin--van Kampen duals are \emph{metrizable}.
An introduction to the pertinent part of this duality theory will be given
in Section 5.

A subset $C$ of a partially ordered set $P$ is \emph{cofinal} (in $P$) if
for each $p\in P$, there is $c\in C$ such that $p\le c$.
In this paper, families of sets are always ordered by $\sub$.

All groups considered in this overview are assumed, without further notice,
to be locally quasiconvex.
This is a mild restriction, meaning that the group admits reasonably many
continuous homomorphisms into the circle group.

A topological space is \IndexPrint{$k_\omega$} if
its topology is determined by a countable cofinal family of compact subsets, i.e.,
there are compact sets $K_1,K_2,\dots\sub X$ such that each compact set $K\sub X$
is contained in some $K_n$, and for each set $U\sub X$ with all $U\cap K_n$ open in $K_n$,
the set $U$ is open in $X$.

Topological abelian groups which are subgroups of the dual of a metrizable groups
are exactly the $k_\omega$ groups.
The class of abelian groups containing open $k_\omega$ subgroups includes,
in addition to all locally compact abelian groups:
\begin{itemize}
\itm[-] all free abelian groups on a compact space, indeed on any $k_\omega$ space;
\itm[-] all dual groups of countable projective limits of metrizable (more generally,
\v{C}ech-complete\footnote{A group $G$ is \emph{\v{C}ech-complete} if it has a compact subgroup
$H$ such that the quotient space $G/H$ is complete metrizable.}) abelian groups;
\itm[-] all dual groups of abelian pro-Lie groups defined by countable systems \cite{glockner_gramlich,hof_mor:book_prolie}.
\end{itemize}
Moreover, this class is preserved by
countable direct sums, closed subgroups, and finite products \cite{glockner_gramlich}.

Consider the set $\NN$ with the partial order $f\le g$ if $f(n)\le g(n)$ for all $n$.
\index{$\NN$, coordinate-wise order $\le$}
The \emph{cofinality} of a partially ordered set $P$, denoted \IndexPrint{$\cof(P)$},
is the minimal cardinality of a cofinal subset of $P$.
The cardinal number \IndexPrint{$\fd$} is the cofinality of $\NN$ with respect to $\le$.
This cardinal was extensively studied \cite{vD, BlassHBK}, and for the present purposes
it may be thought of as some constant cardinal between $\aleph_1$ and the continuum (inclusive).

For a cardinal number $\kappa$, thought of as a set of cardinality $\kappa$,
the set \IndexPrint{$[\kappa]^\alephes$} is the family of all countable subsets of $\kappa$.
The \emph{weight} of a topological space $X$
is the minimal cardinality of a basis of open sets for the topology of $X$.
For brevity, define the \emph{compact weight} of $X$ to be the
supremum of the weights of compact subsets of $X$.
For nondiscrete (locally) compact groups, the character is equal to the (compact) weight.
The main theorem of this paper, stated in an inner language, is the following one.
Note that this theorem is directly applicable to every group containing an open
abelian non-locally compact $k_\omega$ group $G$.

\bthm\label{MainInner}
Let $G$ be an abelian non-locally compact $k_\omega$ group.
Let $\kappa$ be the compact weight of $G$,
and $\lambda$ be the minimum among the compact weights of the quotients of $G$ by compact subgroups.
Then the character of $G$ is the maximum of $\fd$, $\kappa$, and the cofinality of $[\lambda]^\alephes$.
\ethm

In particular, if the group $G$ has no proper compact subgroups (this is the case for
the free abelian groups considered below),
or more generally, if quotients by compact subgroups do not decrease the compact weight of $G$,
then the character of $G$ is the maximum of $\fd$ and $\cof([\kappa]^\alephes)$.

Theorem \ref{MainInner} reduces the computation of the character of the group $G$ to
the purely combinatorial task of estimating the cofinality of $[\lambda]^\alephes$.
The estimation of $\cof([\lambda]^\alephes)$, for a given uncountable cardinal $\lambda$,
is a central goal in Shelah's pcf theory.
The last section of this paper is dedicated to an introduction of this theory and
its applications in our context.
In contrast to cardinal exponentiation, the function $\lambda\mapsto\cof([\lambda]^\alephes)$
is very tame. For example, if there are no large cardinals (in a certain canonical
model of set theory)\footnote{It is not even possible to prove, using the standard axioms of set theory,
that the existence of such cardinals is \emph{consistent}.}, then $\cof([\lambda]^\alephes)$ is simply
$\lambda$ if $\lambda$ has uncountable cofinality, and $\lambda^+$ (the successor of $\lambda$) otherwise.
Thus, the axiom \emph{SSH}, asserting that $\cof([\lambda]^\alephes)\le\lambda^+$, is extremely weak.
Moreover, without any special hypotheses, $\cof([\lambda]^\alephes)$ can be estimated, and in many cases
computed exactly.

For brevity, denote the character of a topological group $G$ by \IndexPrint{$\chi(G)$}.
Following is a summary of consequences of the main theorem.

\bcor\label{manyex}
In the notation of Theorem \ref{MainInner}:
\be
\itm $\chi(G)\le\kappa^\alephes$.
\itm If $\kappa=\kappa^\alephes$, then $\chi(G)=\kappa$.
\itm If $\lambda=\aleph_n$ for some $n$, then $\chi(G)=\max(\fd,\kappa)$.
\itm If $\lambda=\aleph_\mu$, for a limit cardinal $\mu$ below the first fixed point of the $\aleph$ function,
and $\mu$ has uncountable cofinality, then $\chi(G)=\max(\fd,\kappa)$.
\itm If $\lambda=\aleph_\alpha$ is smaller than the first fixed point of the $\aleph$ function, then
$\chi(G)$ is smaller than $\max(\fd^+,\kappa^+,\aleph_{\card{\alpha}^{+4}})$.
\itm If SSH holds, then:
\be
\itm If $\lambda<\kappa$ or $\cof(\lambda)>\alephes$, then $\chi(G)=\max(\fd,\kappa)$.
\itm If $\lambda=\kappa$ and $\cof(\lambda)=\alephes$, then $\chi(G)=\max(\fd,\kappa^+)$.
\ee
\ee
\ecor

The proof of these theorems spans throughout the entire paper, but the paper is designed so that each
reader can read the sections accessible to him or her, taking for granted the other ones.

In Section 2, we set up a general framework for studying bounded sets in topological groups.
The level of generality is just the one needed to capture available methods from the context of
topological vector spaces, and import them to the seemingly different context of separable topological groups
with translations by elements of a dense subset.
This is done in Section 3, which concludes by showing that in metrizable groups, precompact subsets of dense subgroups determine the precompact
subsets of the full group. It follows that the precompact sets in the group and in its
dense subgroup have the same cofinal structure. These are, essentially, the only two results from the first two sections
needed for the remaining sections.
In a first reading of Sections 2 and 3, the reader may wish to consider only the special case of topological
groups with translations by elements of a dense subset, since this is the case needed in the concluding results of these sections.

In Section 4, the approach of Section 3 is generalized from separable to arbitrary metrizable groups.
The \emph{density} of a topological group $G$, \IndexPrint{$\density(G)$}, is the minimal cardinality of a dense subset of that space.
We define the \emph{local density} of $G$, \IndexPrint{$\localdensity(G)$}, to be the minimal density of a neighborhood of the identity
element of $G$. Let \IndexPrint{$\PK(G)$} denote the family of all precompact subsets of $G$.
The main result of this section is the following theorem.
In this theorem, which is of independent interest, we do not require that $G$ is
locally quasiconvex or abelian.

\bthm\label{goodnewsthm}
Let $G$ be a metrizable non-locally precompact group.
The cofinality of $\PK(G)$ is equal to the maximum of $\fd$, $\density(G)$, and $\cof([\localdensity(G)]^\alephes)$.
\ethm

In Section 5 we use Theorem \ref{goodnewsthm} and methods of Pontryagin--van Kampen duality to
prove the following theorem. A topological abelian group is \emph{complete} if it is
complete with respect to its uniformity. (Being abelian, the left, right, and two-sided
uniformities of the group coincide.)

\bthm\label{MainOuter}
Let $G$ be a complete abelian group whose dual group is a metrizable non-locally precompact group $\Gamma$.
Then $\chi(G)$ is the maximum of $\fd$, $\density(\Gamma)$, and $\cof([\localdensity(\Gamma)]^\alephes)$.
\ethm

This already puts us in a position to prove, in Section 6, the following result.
We state it in full because the estimations are slightly simpler than those
in Corollary \ref{manyex}.

\bthm\label{thm:MchiAX1}
Let $X$ be a nondiscrete $k_\omega$ space.
Let $\kappa$ be the compact weight of $X$. Then the character of $A(X)$ is the maximum of $\fd$
and $\cof([\kappa]^\alephes)$.
\ethm

\bcor\label{cor:MchiAX2}
In the notation of Theorem \ref{thm:MchiAX1}:
\be
\itm $\chi(A(X))\le\kappa^\alephes$, and if $\kappa=\kappa^\alephes$, then $\chi(A(X))=\kappa$.
\itm If $\kappa=\aleph_n$ for some $n\in\N$, then $\chi(A(X))=\max(\fd,\aleph_n)$.
\itm If $\kappa=\aleph_\mu$, for $\mu$ smaller than the first fixed point of the $\aleph$ function,
and $\mu$ is a limit cardinal of uncountable cofinality, then $\chi(A(X))=\max(\fd,\aleph_\mu)$
\itm If $\kappa=\aleph_\alpha$ is smaller than the first fixed point of the $\aleph$ function, then
$\chi(A(X))$ is smaller than $\max(\fd^+,\aleph_{\card{\alpha}^{+4}})$.
\itm If SSH holds, then:
\be
\itm If $\cof(\kappa)>\alephes$, then $\chi(A(X))=\max(\fd,\kappa)$.
\itm If $\cof(\kappa)=\alephes$, then $\chi(A(X))=\max(\fd,\kappa^+)$.
\ee
\ee
\ecor

By virtue of \cite[Corollary 2.3]{NT2}, Theorem \ref{thm:MchiAX1} also holds
for the free \emph{nonabelian} topological group $F(X)$.

The result in Theorem \ref{thm:MchiAX1} was previously known only in few of the cases covered by this theorem 
\cite{NT1,NT2}, for example when $X$ is compact,
or when, in addition to the premise in our theorem, all compact subsets of $X$ are metrizable
\cite{NT2}.
However, Theorem \ref{thm:MchiAX1} does not capture all of the related results of
\cite{NT1,NT2}. The proofs in \cite{NT1,NT2} are more combinatorially oriented than ours.

In Section 7 we develop the remaining Pontryagin--van Kampen theory required to
deduce Theorem \ref{MainInner} from Theorem \ref{MainOuter}.
Section 8 introduces and applies pcf theory, to obtain the concrete estimations in Corollary
\ref{manyex} and Corollary \ref{cor:MchiAX2}.

We note that all estimations in Corollary \ref{manyex} apply to Theorem \ref{MainOuter} as well,
which may be viewed by some readers as the main result of this paper.

\section{Bounded sets in topological groups}

The unifying concept of this paper is that of boundedness in topological groups.
This concept plays a central role in a number of studies in functional analysis and topology.
In its most abstracted form, a \emph{boundedness} (or \emph{bornology}) on a topological space
$X$ is a family of subsets of $X$ that is closed under taking subsets and unions of finitely many elements, and
contains all finite subsets of $X$.\footnote{In set theoretic terms, this defines a (not necessarily proper)
\emph{ideal} on $X$ containing all singletons.}
The abstract approach has found applications in several areas of mathematics -- see
the introduction and references in \cite{Beer08}. In particular, Vilenkin \cite{vilenkin}
applied this approach in the realm of topological groups.
Here, we focus on well-behaved boundedness notions in topological groups, which make
it possible to simultaneously extend some earlier studies in locally convex topological vector spaces
as well as seemingly unrelated studies of general topological groups.

We use the following notational conventions throughout the paper.
For a set $X$, let \IndexPrint{$P(X)$} denote the family of all subsets of $X$,
and let \IndexPrint{$\Fin{X}$} denote the family of all \emph{finite} subsets of $X$.
An \emph{operator} $t$ on $P(X)$ is a function $t\colon P(X)\to P(X)$.
Throughout, $G$ is an infinite Hausdorff topological group with identity element $e$ (or $0$ if $G$ is restricted to be abelian),
and $T$ is a set of operators on $P(G)$.

\bdfn
For an operator $t$ on $P(G)$ and $A\sub G$, write \IndexPrint{$t\mult A$} for $t(A)$.
Let $T$ be a set of operators on $P(G)$.
\be
\itm For $H\sub T$, let \IndexPrint{$H\mult A:=\Union_{t\in H}t \mult A$}.
\itm A set $B\sub G$ is \emph{$T$-bounded} (\emph{bounded}, when $T$ is clear from the context)
if for each neighborhood $U$ of $e$ there is a finite set $F\sub T$ such that $B\sub F\mult U$.
\ee
\edfn

The following axioms guarantee that the family of $T$-bounded sets is a boundedness notion.

\bdfn\label{bmap}
A \emph{boundedness system} is a pair $(G,T)$
such that $G$ is a topological group, $T$ is a set of operators on $P(G)$,
and the following axioms hold:
\be
\itm[(B1)] For each open set $U$ and each element $t\in T$, the set $t\mult U$ is open;
\itm[(B2)] For each neighborhood $U$ of $e$, we have that $T\mult U=G$;
\itm[(B3)] For each $T$-bounded set $A\sub G$ and each $t\in T$, the set $t\mult A$ is $T$-bounded;
\itm[(B4)] For all $A\sub B\sub G$ and each $t\in T$, we have that $t\mult A\sub t\mult B$;
\itm[(B5)] For each $S\sub T$ with $\card{S}<\card{T}$, there is a neighborhood $U$ of $e$ such that $S\mult U\neq G$;
\itm[(B6)] For each $n$, there is a neighborhood $U$ of $e$ such that for all $F\sub T$ with $\card{F}\le n$,
we have that $F\mult U\neq G$.
\ee
A boundedness system $(G,T)$ is said to be \emph{metrizable} if $G$ is metrizable.
\edfn

Axiom (B5) is assumed since one can restrict attention to
a set $T'\sub T$ of minimal cardinality such that $T'\mult U=G$ for each neighborhood $U$ of $e$.
Axiom (B6) is added to avoid trivialities. By moving to the semigroup of operators generated by $T$,
we may assume that $T$ is a semigroup. We will, however, not make use of this fact.

The following example shows that precompact sets need not be bounded when $G$ is not complete.
However, we have the subsequent Lemma \ref{hh}.
\bexm
Consider the additive group $\bbQ$ of rational numbers, equipped with its standard topology.
Enumerate $\bbQ$ as $\set{q_n}{n\in\N}$
and let $\{x_n\}$ be a sequence of rational numbers converging to $\sqrt{2}$.
Taking $T=\N$, we define
$n\mult A= (q_n+A)\sm\set{x_k}{k\ge n}$.
Then the sequence $\{x_n\}$ is an precompact but unbounded subset of $\bbQ$.
\eexm

\blem\label{hh}
For each boundedness system $(G,T)$:
\be
\itm Every compact set $K\sub G$ is bounded.
\itm The family of bounded subsets of $G$ is a boundedness.\qed
\ee
\elem

The following two examples of boundedness systems are well known.
In these examples, we identify $T$ with some set of parameters defining
the elements of $T$. In general, we may identify $T$ with any set
$S$ of the same cardinality, by modifying the definition of $\mult$ appropriately.

\bexm[\IndexPrint{Standard boundedness on topological vector spaces}]\label{An}
Let $E$ be a topological vector space. Take $T=\N$, and
define $n\mult A=\set{nv}{v\in A}$ for each $A\sub V$.
For example, Axiom (B2) holds since $\lim_n \frac{1}{n}v=\vec{0}$ for
each $v\in E$. The $\N$-bounded sets are those bounded in the ordinary sense.
\eexm

In Example \ref{An}, if $E$ is a locally convex topological vector space,
we may alternatively define $n\mult A=nA=\set{v_1+\dots+v_n}{v_1,\dots,v_n\in A}$ for each $A\sub V$,
and obtain the same bounded sets.
More generally, for any connected multiplicative topological group $G$,
we can take $T=\N$ and $n\mult A=A^n=\set{a_1a_2\cdots a_n}{a_1,a_2,\dots,a_n\in A}$.
Let $U$ be an open and symmetric neighborhood of $e$.
Then $\N\mult U$ is an open, and therefore also closed, subgroup of $G$.
Thus, $\N\mult U=G$.

\bexm[\IndexPrint{Standard boundedness on Topological groups}]\label{nbg}
Fix any dense subset $T$ of $G$ of minimal cardinality.
Define $t\mult A=tA=\set{ta}{a\in A}$ for all $t\in T, A\sub G$.
The $T$-bounded sets are the precompact subsets of $G$. Axiom (B6) holds because our groups are
assumed to be infinite Hausdorff.
Indeed, let $x_1,\dots, x_{n+1}$ be distinct elements of $G$.
Take a symmetric neighborhood $U$ of the identity element
such that $x_iU^2\cap x_jU^2=\emptyset$ for all distinct $i$ and $j$.
Assume that $F\sub G$, $\card{F}\leq n$ and $FU=G$.
Then there are an element $a\in F$ and distinct indices $i$ and $j$
such that $\{x_i,x_j\}\subseteq aU$.
Then $x_j\in x_iU^2$; a contradiction.
Axiom (B2) is equivalent to the density of $T$:
If $U$ is a symmetric neighborhood of the identity element, then
$t\in T\cap (gU)$ if and only if $g\in tU$.
The remaining axioms are a straightforward consequences
of basic properties of topological groups.

It follows that if $T\sub G$ is a set of translations
then $(G,T)$ is a boundedness system if and only if $T$ is dense in $G$.
\eexm

When a topological group also happens to be a topological vector space, the term
\emph{standard boundedness system on $G$} has two contradictory interpretations.
When we wish to use the one of topological vector spaces, we will say so explicitly.

The two canonical examples were combined by Hejcman \cite{Hejcman},
who considered the case $T=D\x\N$, where $D$ is a dense subset of $G$, and $(d,n)\mult A=dA^n$.
The $T$-bounded sets are  the standard bounded sets when $G$ is a topological vector space, and the precompact
sets when $G$ is a locally compact group.

\bdfn\label{kb}
Let $(G,T)$ be a boundedness system and $\kappa$ be an infinite cardinal number.
A set $A\sub G$ is \emph{$\kappa$-bounded}
(with respect to $T$) if, for each neighborhood $U$ of $e$, there is
a set $S\sub T$ of cardinality at most $\kappa$ such that $A\sub S\mult U$.
The \emph{boundedness number} of $A$ in $(G,T)$, denoted \IndexPrint{$\boundedness_T(A)$},
is the minimal cardinal $\kappa$ such that $A$ is $\kappa$-bounded.
\edfn

Axiom (B6) asserts that $\boundedness_T(G)\ge\alephes$.

\bdfn
For a topological group $G$ and a set $A\sub G$, \IndexPrint{$\boundedness(A)$} is the minimal cardinal $\kappa$ such that
for each neighborhood $U$ of $e$, there is $S\sub A$ such that $\card{S}\le\kappa$, and $A\sub SU$.
\edfn

For the standard boundedness system $(G,T)$ on a topological group $G$ (Example \ref{nbg}),
the cardinal $\boundedness_T(G)$ does not depend on the choice of the dense subset $T$. Indeed, we have the following.

\blem[folklore]\label{subgp-kbdd}
Let $(G,T)$ be a standard boundedness system on $G$. Then:
\be
\itm $\boundedness_T(A)=\boundedness(A)$ for all $A\sub G$.
\itm If $A\sub B\sub G$, then $\boundedness(A)\le \boundedness(B)$.
\ee
\elem
\bpf
(2) Clearly, $\boundedness_T(A)\le \boundedness_T(B)$. Thus, it suffices to prove (1).

$(\ge)$ Fix a neighborhood $U$ of $e$ in $G$.
Let $V$ be a neighborhood of $e$ in $G$, such that $V=V\inv$ and $V^2\sub U$.
Let $S\sub T$ be such that $\card{S}\le \boundedness_T(A)$, and $A\sub SV$.
By thinning out $S$ if needed, we may assume that for each $s\in S$, the set $sV$ intersects $A$.
For each $s\in S$, pick an element $a_s\in sV\cap A$. Then $s\in a_sV$, and thus $sV\sub a_sV^2\sub a_sU$.
Let $S'=\set{a_s}{s\in S}$.
Then $S'\sub A$, $\card{S'}\le\card{S}\le \boundedness_T(A)$, and $A\sub SV\sub S'U$.

$(\le)$ Similar, using that $T$ is dense in $G$.
\epf

\bcor
For a standard boundedness system $(G,T)$ on a topological group, the cardinality of $T$ is $\density(G)$.\qed
\ecor

Thus, if $(G,T)$ is a boundedness system with $G$ a $\sigma$-compact group, then $\boundedness_T(G)=\alephes$.
But if $G$ is (nonmetrizable and) not separable, then for the
standard boundedness system on $G$, $\card{T}=\density(G)>\alephes$. That is, for each neighborhood $U$ of $e$ there
is a countable $S\sub T$ such that $S\mult U=G$, but there is no such $S$ independent on $U$.

Recall that for infinite cardinals $\kappa$ and $\lambda$, $\kappa\cdot\lambda=\max(\kappa,\lambda)$.

\bprp\label{bsandwich}
Let $(G,T)$ be a boundedness system. Then
$$\boundedness_T(G)\le \card{T}\le\chi(G)\cdot \boundedness_T(G).$$
In particular:
\be
\itm For metrizable $G$, $\card{T}=\boundedness_T(G)$.
\itm $\boundedness(G)\le \density(G)\le\chi(G)\cdot \boundedness(G)$.
\itm For metrizable $G$, $\boundedness(G)=\density(G)$.
\ee
\eprp
\bpf
$\card{T}\le\chi(G)\cdot \boundedness_T(G)$:
Let $\set{U_\alpha}{\alpha<\chi(G)}$ be a neighborhood base of $G$ at $e$.
For each $\alpha<\chi(G)$, let $S_\alpha\sub T$ be such that $\card{S_\alpha}\le \boundedness_T(G)$,
and $S_\alpha\mult U_\alpha=G$. Let $S=\Union_{\alpha<\chi(G)}S_\alpha$.
For each neighborhood $U$ of $e$, $S\mult U=G$. It follows that
$\card{T}=\card{S}\le \chi(G)\cdot \boundedness_T(G)$.

For (2) and (3), consider the standard boundedness system on $G$.
\epf

Thus, when considering metrizable groups, we may replace $\boundedness_T(G)$ by $\card{T}$, or by $\density(G)$
when the standard boundedness system is considered.

We give some examples, using the multiplicative \emph{torus group} $\bbT=\set{z\in\bbC}{\card{z}=1}$\index{$\bbT$}.

\bexm
The inequalities in Proposition \ref{bsandwich} cannot be improved, not even for the standard
boundedness system (Proposition \ref{bsandwich}(3)) on powers of the torus:
For compact groups $G$ of cardinality $2^\kappa$, we have that $\boundedness(G)=\alephes$,
and $\density(G)=\log(\kappa)$, where
\IndexPrint{$\log(\kappa)$} is defined as $\min\set{\lambda}{\kappa\le 2^\lambda}$ \cite[Theorem 3.1]{ComfortHBK}.

Thus, for an infinite cardinal $\kappa$, we have that
$\boundedness(\bbT^\kappa)=\alephes$, $\density(\bbT^\kappa)=\log(\kappa)$
and $\chi(\bbT^\kappa)=\kappa$. The inequality $\alephes\le\log(\kappa)\le\kappa$ cannot be
improved. Indeed, for $\fc:=2^\alephes$, we have the following:
\be
\itm $\kappa=\alephes$ gives $\boundedness(G)=\density(G)=\chi(G)=\alephes$.
\itm $\kappa=\fc$ gives $\boundedness(G)=\density(G)=\alephes<\chi(G)=\fc$.
\itm $\kappa=\fc^+$ gives $\boundedness(G)=\alephes<\density(G)=\log(\fc^+)<\chi(G)=\fc^+$.
\itm $\kappa=\beth_\w$ gives $\boundedness(G)=\alephes<\density(G)=\chi(G)=\beth_\w$.
\ee
Here, the cardinal \IndexPrint{$\beth_\w$} is defined as the supremum of all cardinals \IndexPrint{$\beth_n$}, $n\in\N$,
where $\beth_1=2^{\alephes}$ and for each $n>1$, $\beth_n=2^{\beth_{n-1}}$.
\eexm

\section{When $T$ is countable}

Methods and ideas from the context of topological vector spaces,
developed by Saxon and S\'anchez--Ruiz \cite{SaxSan95}, and by Burke and Todorcevic \cite{BurTod96},
generalize to general boundedness systems $(G,T)$ with $T$ countable.
Even for the standard boundedness systems on topological groups,
some of the obtained results were apparently not observed earlier.

\bdfn
A boundedness system $(G,T)$ is \emph{locally bounded}
if there is in $G$ a neighborhood base at $e$, consisting of bounded
sets.
\edfn


Let $P$ and $Q$ be partially ordered sets. We write \IndexPrint{$P\preceq Q$} if there is an order
preserving $f\colon P\to Q$ with image cofinal in $Q$.
We say that $P$ is \emph{cofinally equivalent} to $Q$ if $P\preceq Q$ and $Q\preceq P$.
Our notion of cofinal equivalence is stronger and simpler than the standard one.
This variation will not affect our results.

If $P\preceq Q$, then $\cof(Q)\le\cof(P)$.

\bdfn
Let $(G,T)$ be a boundedness system.
\IndexPrint{$\Bdd_T(G)$} is the family of $T$-bounded subsets of $G$.
$\Bdd_T(G)$ is partially ordered by the relation $\sub$.
When $(G,T)$ is a standard boundedness system,
$\Bdd_T(G)$ is the family of precompact subsets of $G$, which we denote
for simplicity by \IndexPrint{$\PK(G)$}.
\edfn

\brem
If $G$ is $T$-bounded, then $\Bdd_T(G)$ is cofinally equivalent to the singleton $\{1\}$.
\erem

For a function $f\colon X\to Y$ and sets $A\sub X$ and $B\sub Y$, we use the notation
$f[A]=\set{f(a)}{a\in A}$ and $f\inv[B]=\set{x\in X}{f(x)\in B}$.
\index{$f[A]$}\index{$f\inv[A]$}

For locally convex topological vector spaces with the standard boundedness structure,
the following is pointed out in \cite[Theorem 2.5]{BurTod96}.
Recall that when $T$ is countable, we may identify $T$ with $\N$.

\bprp\label{cofeqN}
If a boundedness system $(G,\N)$ is locally bounded and $G$ is unbounded, then $\Bdd_\N(G)$ is cofinally equivalent to $\N$.
\eprp
\bpf
Fix a bounded neighborhood $U$ of $e$, such that for each finite $F\sub\N$, $F\mult U\neq G$.
Define $\vphi\colon G\to\N$ by
$$\vphi(g)=\min\set{n}{g\in n\mult U}.$$
The functions $K\mapsto\max\vphi[K]$ and $n\mapsto \vphi\inv[\{1,\dots,n\}]$
establish the required cofinal equivalence.
\epf

Systems which are \emph{not} locally bounded are more interesting in this respect.
Assume that $(G,\mathbb N)$ is a metrizable boundedness system,
and let $U_n$, $n\in\N$, be a neighborhood base at $e$.

\bdfn\label{Psidef}
Define a map $\Psi\colon G\to\NN$\index{$\Psi$} by
$$x \mapsto \vphi_x(n) = \min\set{m}{x\in m\mult U_n}.$$
\edfn
For a bounded set $B\sub\NN$, the function $f:=\max B\in\NN$ is defined by
$f(n)=\max\set{g(n)}{g\in B}$.
Define functions
$\Bdd_\N(G)\to\NN$ and $\NN\to\Bdd_\N(G)$, respectively, by
\beq
K & \mapsto & \max\Psi[K];\\
f & \mapsto & \Psi\inv[\set{g\in\NN}{g\le f}].
\eeq
Both functions are monotone, and the image of the latter is cofinal in $\Bdd_\N(G)$.

For locally convex topological vector spaces with the standard boundedness structure, the following is proved
in \cite[Proposition 1]{SaxSan95} and in \cite[Theorem 2.5]{BurTod96}.

\bthm\label{cofeqNN}
Let $(G,\N)$ be a metrizable non-locally bounded boundedness system.
Then $\Bdd_\N(G)$ is cofinally equivalent to $\NN$.
\ethm
\bpf
As compact sets are bounded, it suffices to show that there is a neighborhood base $U_n$, $n\in\N$, at $e$,
and for each $f\in\NN$, there is a compact set $K\sub G$ such that $f\le\max\Psi[K]$.

Let $U_n$, $n\in\N$, be a descending neighborhood base at $e$. As $U_1$ is not bounded, we may assume
(by shrinking $U_2$ if needed) that there is no $m$ such that $U_1\sub\{1,\dots,m\}\mult U_2$.
Continuing in the same manner, we may assume that for each $n$, there is no $m$ such that
$U_n\sub\{1,\dots,m\}\mult U_{n+1}$.

Given $f\in\NN$, choose for each $n$ an element $x_n\in U_n\sm \{1,\dots,f(n)\}\mult U_{n+1}$.
As the original sequence $U_n$ was descending to $e$, the elements $x_n$ converge to $e$, and thus
the set $\set{x_n}{n\in\N}\cup\{e\}$ is compact, as required.
\epf

\bcor\label{Polish}
Let $G$ be a separable metrizable non-locally precompact group.
Then $\PK(G)$ is cofinally equivalent to $\NN$.\qed
\ecor

\bdfn
For a topological space $X$, let \IndexPrint{$C(X,\bbT)$} be the group of all continuous functions
from $X$ into $\bbT$, with pointwise multiplication, endowed with the \emph{compact-open topology}.
That is, a neighborhood base at the constant function $1$ is given by the sets
$$\set{f\in C(X,\bbT)}{\card{f(x)-1}<\epsilon\mbox{ for all }x\in K},$$
where $K$ is a compact subset of $X$, and $\epsilon$ is a positive real number.
\edfn

A \emph{Polish group} is a complete, separable, metrizable group. We give
two well known examples of non-locally compact Polish groups, and
where it is not immediately clear (without Corollary \ref{Polish}) that $\PK(G)$ is cofinally equivalent to $\NN$.

\bexm\label{exCKT}
Let $L$ be a Lie group, for example $\bbT$ or the group of unitary $n\x n$ complex matrices.
Let $K$ be a compact metric space.
The group $C(K,L)$ is Polish, with the metric given by the supremum norm.
$C(K,L)$ is not locally compact (unless $K$ is finite).
By Theorem \ref{cofeqNN}, the family of compact subsets
of $C(K,L)$ is cofinally equivalent to $\NN$.
\eexm

\bexm\label{exSN}
Consider the group \IndexPrint{$S_\N$} of permutations on $\N$, where for each finite
$F\sub\N$, the set $U_F$ of all permutations fixing $F$ is a neighborhood of the identity.
This defines a neighborhood base at the identity permutation, and thus a topology on $S_\N$.
The nonabelian group $S_\N$ is Polish and non-locally compact. Thus, its compact subsets
are cofinally equivalent to $\NN$.
\eexm


For functions $f,g\in\NN$, the notation $f\le^* g$\index{$\le^*$} stands for $f(n)\le g(n)$ for all but finitely many $n$.
The cardinal number \IndexPrint{$\fb$} is the minimal cardinality of a $\le^*$-unbounded subset of $\NN$.
The cardinal $\fb$ is uncountable, and can consistently be any regular uncountable cardinal
not larger than $\mathfrak c$. More details about this cardinal
are available in \cite{vD, BlassHBK}.

For locally convex topological vector spaces  with the standard boundedness structure, the following is Corollary 2.6 of \cite{BurTod96}.

\bcor\label{sigmabdd}
Let $(G,\N)$ be a metrizable boundedness system.
\be
\itm For each family $\cF\sub\Bdd_\N(G)$ with $\card{\cF}<\fb$, there is a countable family $\cS\sub\Bdd_\N(G)$
such that each member of $\cF$ is contained in a member of $\cS$.
\itm Each union of less than $\fb$ bounded subsets of $G$ is a union of countably many bounded subsets of $G$.
\ee
\ecor
\bpf
The assertions are immediate when $G$ is locally bounded. Thus, assume it is not.
Then (1) follows from the cofinal equivalence of $\Bdd_\N(G)$ and $\NN$, and (2) follows from (1).
\epf

\bdfn
A group $G$ is \emph{metrizable modulo precompact}
if there is a precompact subgroup $K$ of $G$,
such that the coset space $G/K$ is metrizable.
\edfn

\bexm
All \v{C}ech-complete groups, and all almost-metrizable groups,
are metrizable modulo precompact.
\eexm

For a nonabelian group $G$, the coset space $G/K$ need not be a group since we do not require $K$ to be a \emph{normal} subgroup.
However, the concept of boundedness extends naturally to the coset space $G/K$, and we have the following.

\blem\label{modprec}
Let $K$ be a precompact subgroup of $G$, and $\pi\colon G\to G/K$ be the canonical quotient map.
\be
\itm If $P\in\PK(G)$, then $\pi[P]\in\PK(G/K)$.
\itm If $Q\in\PK(G/K)$, then $\pi\inv[Q]\in\PK(G)$.
\itm $\PK(G)$ is cofinally equivalent to $\PK(G/K)$.
\ee
\elem
\bpf
(1) Precompactness of $K$ is not needed here:
Let $U$ be a neighborhood of $eK$ in $G/K$.
As $\pi\inv[U]$ is a neighborhood of $e$ in $G$, there is a finite $F\sub G$ such that
$P\sub F\pi\inv[U]$. Then $\pi[P]\sub\pi[F\pi\inv[U]]=FU$.

(2) Let $U$ be a neighborhood of $e$ in $G$. Take a neighborhood $W$ of $e$ such that $W^2\sub U$.
As $K$ is precompact, there is a neighborhood $V$ of $e$ such that $VK\sub KW$.\footnote{This is standard:
Take a neighborhood $W_0$ of $e$ with $W_0^2\sub W$, and then
take a finite $F\sub K$ such that $K\sub FW_0$. For each $g\in F$, $e\cdot g=g\in FW_0$, and thus
there is a neighborhood $V_g$ of $e$ with $V_g\cdot g\sub FW_0$. Take $V=\bigcap_{g\in F}V_g$. Then $VF\sub FW_0$, and thus
$VK\sub VFW_0\sub FW_0W_0\sub FW$.}
As $K$ is precompact, there is a finite $I\sub G$ such that $K\sub IW$.

The set $\pi[V]$ is a neighborhood of $eK$ in $G/K$.
Take a finite subset $F$ of $G$ such that $Q\sub \pi[F]\pi[V]$.
Then $\pi\inv[Q]\sub \pi\inv[\pi[F]\pi[V]]=FKVK\sub FK^2W=FKW\sub FIW^2\sub FIU$, and $FI$ is finite.

(3) If $P\in\PK(G)$, then $Q=\pi[P]\in\PK(G/K)$, and $\pi\inv[Q]\in\PK(G)$, and contains $P$.
Thus, the map $Q\mapsto\pi\inv[Q]$ shows that $\PK(G/K)\preceq \PK(G)$.
Similarly, if $Q\in\PK(G/K)$, then $P=\pi\inv[Q]\allowbreak\in\PK(G)$, and $Q=\pi[P]\in\PK(G/K)$,
and thus the map $P\mapsto\pi[P]$ gives $\PK(G)\preceq\PK(G/K)$.
\epf

\bcor\label{Baire}
Let $G$ be a separable, metrizable modulo precompact, Baire group.
If $G$ is a union of fewer than $\fb$ precompact sets, then $G$ is locally precompact.
\ecor
\bpf
By Lemma \ref{modprec}, we may assume that $G$ is metrizable.
By Corollary \ref{sigmabdd}, $G$ is a union of countably many precompact sets.
As the closure of precompact sets is precompact, we may assume that these sets are closed.
As $G$ is Baire, one of these sets has nonempty interior. It follows that there is a precompact
neighborhood of $e$.
\epf

If every bounded subset of a normed space is separable, then the space is separable.
Dieudonn\'e \cite{D} asked to what extent this can be generalized to locally convex
topological vector spaces. Burke and Todorcevic answered this question completely, by
showing that the same assertion holds in all locally convex
topological vector spaces if, and only if, $\aleph_1<\fb$ \cite{BurTod96}.
One direction of this assertion is generalized by the following theorem.
This theorem, which is trivial when applied to standard boundedness systems on topological
groups, is nontrivial in general.

\bthm\label{bddsepallsep}
Let $(G,\N)$ be a metrizable boundedness system with $\density(G)<\fb$.
If all bounded subsets of $G$ are separable, then $G$ is separable.
\ethm
\bpf
Assume otherwise, and let $D$ be a discrete subset of $G$ of cardinality $\aleph_1$.
As $\aleph_1<\fb$, we have by Corollary \ref{sigmabdd} that $D$ is a union of countably many bounded sets.
Thus, $D$ has a (discrete, of course) bounded subset of cardinality $\aleph_1$.
\epf

\bprp\label{subseq}
For each sequence $x_n\to x$ in $G$, there is a subsequence $\{y_n\}$ of $\{x_n\}$ such that
$\vphi_{y_n}$ converges to a function $f\le \vphi_x$.
\eprp
\bpf
For each $k$, the set $\set{y\in G}{\varphi_y(k) \le \varphi_x(k)}$ is
an open neighborhood of $x$. 
Thus, $\vphi_{x_n}(1)\le \vphi_x(1)$ for all but finitely many $n$.
Therefore, there is $m_1\le \vphi_x(1)$ such that $I_1=\set{n}{\vphi_{x_n}(1)=m_1}$ is infinite.

Inductively, given the infinite $I_{k-1}\sub\N$, we have
that $\vphi_{x_n}(k)\le \vphi_x(k)$ for all but finitely many $n\in I_{k-1}$, and thus
there is $m_k\le \vphi_x(k)$ such that $I_k=\set{n\in I_{k-1}}{\vphi_{x_n}(k)=m_k}$ is infinite.

For each $k$, pick $i_k\in I_k$ with $i_k>i_{k-1}$. Then $\vphi_{x_{i_k}}\to f$, where $f(k)=m_k\le \vphi_x(k)$ for all $k$.
\epf

The next result tells that if the group has a small dense subset, then
the bounded subsets of its completion are determined by the bounded subsets
of any dense subgroup of $G$.
A special case of it was proved by Grothendieck \cite{G}, and extended in
\cite[Theorem 2.1]{BurTod96}, for $G$ a separable metrizable locally convex topological vector space.

\bthm\label{catchcompact}
Let $(G,\N)$ be a metrizable boundedness system with $\density(G)<\fb$.
Let $D$ be a dense subset of $G$. For each bounded $K\sub G$, there is
a bounded $J\sub D$ such that $K\sub \cl{J}$.
\ethm
\bpf
Assume that $G$ is locally compact, and let $U$ be a compact neighborhood of $e$.
Take a finite $F\sub\N$ such that $K\sub F\mult U$, and let $J=D\cap (F\mult U)$.
Then $K\sub \cl{J}$.

Next, assume that $G$ is not locally compact.
As $\density(G)<\fb$, there is $K'\sub K$ such that $\card{K'}<\fb$ and $K\sub\cl{K'}$.
For each $x\in K'$, let $\{x_n\}$ be a sequence in $D$ converging to
$x$. By Proposition \ref{subseq}, we may assume that $\{\vphi_{x_n}\}$
converges to a function $\vphi_x'\le \vphi_x$. The set $\set{x_n}{n\in\N}\cup\{x\}$ is compact,
and thus bounded. Take $g_x$ such that $\vphi_{x_n}\le g_x$ for all $n$.

As $\card{K'}<\fb$, there is $h\in\NN$ such that $g_x\le^* h$ for all $x\in K'$.
We require also that all elements of $\Psi[K]$ are $\le h$.
For each $x\in K'$, we have that $\vphi_{x_n}\le h$ for all but finitely many $n$.
Indeed, let $N$ be such that $g_x(k)\le h(k)$ for all $k>N$.
For all but finitely many $n$,
$$\vphi_{x_n}(1)=\vphi_x'(1)\le \vphi_x(1)\le h(1), \dots, \vphi_{x_n}(N)=\vphi_x'(N)\le \vphi_x(N)\le h(N),$$
as $x\in K$, and for $k>N$, $\vphi_{x_n}(k)\le g_x(k)\le h(k)$.
Thus, for $J=D\cap\Psi\inv[\set{f\in\NN}{f\le h}]$, we have that $K'\sub\cl{J}$,
and therefore also $K\sub\cl{J}$.
\epf

It seems that the following special case of Theorem \ref{catchcompact} was not noticed before.

\bcor\label{groupcathcompact}
Let $G$ be a metrizable group with a dense subgroup $H$.
For each precompact set $K\sub G$, there is a precompact set $J\sub H$ such that $K\sub \cl{J}$.
\ecor
\bpf
As $K$ is precompact and $G$ is metrizable, $K$ is separable. As $H$ is dense in $G$
and $K$ is separable, there is a countable $D\sub H$ such that $K\sub\cl{D}$.
We may assume that $D$ is a group.
Let $G'=\cl{D}$, and apply Theorem \ref{catchcompact}
to $G'$ and $D$ to obtain a bounded set $J\sub D$ such that $K\sub \cl{J}$.
\epf

\bexm
Consider the permutation group $S_\N$ from Example \ref{exSN}.
By Corollary \ref{groupcathcompact}, each compact subset of $S_\N$
is contained in the closure of some precompact set of finitely supported permutations.
\eexm

\brem
There is no assumption on the density of $G$ in corollary \ref{groupcathcompact}.
However, metrizability is needed: A \emph{$P$-group} is a group where every $G_\delta$ set
is open. For each $P$-group $G$ with a proper dense
subgroup $H$, and each $g\in G$, the singleton $\{g\}$ is not contained in the closure of any
precompact subset of $H$. Indeed, if $B\sub H$ is precompact, then $\cl{B}$ is
a compact subset of $G$, and thus finite (countably infinite subsets of $P$-spaces
are closed and discrete), and thus $\cl{B}\sub H$.

For a concrete example, let $\Z_2$ be the two element group, and take $G=(\Z_2)^\kappa$ for some
$\kappa>\alephes$, with the countable box topology, and let $H$ be the group of all $g\in(\Z_2)^{\kappa}$
which are supported on a countable set.
\erem

Corollary \ref{groupcathcompact} implies the following.

\bcor\label{completion}
Let $G$ be a metrizable group with a dense subgroup $H$.
Then $\PK(H)$ is cofinally equivalent to $\PK(G)$.\qed
\ecor

\section{The cofinality of the family of bounded sets}

For locally convex topological vector spaces with the standard boundedness structure,
the following corollary is proved in \cite[Theorem 1]{SaxSan95} and in \cite[Theorem 2.5]{BurTod96}.
In its general form, it follows from Proposition \ref{cofeqN} and Theorem \ref{cofeqNN}.

\bcor\label{hi}
Let $(G,\N)$ be a boundedness system.
\be
\itm If $G$ is bounded, then $\cof(\Bdd_\N(G))=1$.
\itm If $G$ is locally bounded and unbounded, then $\cof(\Bdd_\N(G))=\alephes$.
\itm If $G$ is metrizable non-locally bounded, then $\cof(\Bdd_\N(G))=\fd$.\qed
\ee
\ecor

\blem\label{kappalecofbdd}
Let $(G,T)$ be a boundedness system.
\be
\itm If $G$ is bounded, then $\cof(\Bdd_T(G))=1$.
\itm If $G$ is unbounded, then:
\be
\itm $\alephes\le\cof(\Bdd_T(G))$.
\itm $\boundedness_T(G)\le\cof(\Bdd_T(G))$.
\itm If $\chi(G)\le \card{T}$ (in particular, for metrizable $G$), then $\card{T}\le\cof(\Bdd_T(G))$.
\ee
\ee
\elem
\bpf[Proof of (2)]
(a) Otherwise, $G$ is the union of finitely many bounded sets, and is thus bounded.

(b) Let $\kappa=\cof(\Bdd_T(G))$. By (a), $\kappa\ge\alephes$.
Let $\set{K_\alpha}{\alpha<\kappa}$ be cofinal in $\Bdd_T(G)$.
For each neighborhood $U$ of $e$, there are finite $F_\alpha\sub T$, for $\alpha<\kappa$,
such that $K_\alpha\sub F_\alpha\mult U$. Let $S=\Union_{\alpha<\kappa}F_\alpha$.
Then $\card{S}=\kappa$, and the set $S\mult U$ contains the set $\Union_{\alpha<\kappa}K_\alpha=G$.

(c) Apply (b) and Proposition \ref{bsandwich}.
\epf

\blem\label{lb}
\mbox{}
\be
\itm Let $(G,T)$ be an unbounded locally bounded metrizable boundedness system. Then $\cof(\Bdd_T\allowbreak(G))=\card{T}$.
\itm For each metrizable nonprecompact locally precompact group $G$, we have that $\cof(\PK(G))=\density(G)$.
\ee
\elem
\bpf[Proof of (1)]
Let $U$ be a bounded neighborhood of $e$. Then the set $\set{F\mult U}{F\in \Fin{T}}$ is cofinal
in $\Bdd_T(G)$, and thus $\cof(\Bdd_T(G))\le\card{\Fin{T}}=\card{T}$.
Apply Lemma \ref{kappalecofbdd}.
\epf

\bdfn
For a set $X$,
\IndexPrint{$\Fin{X}^\N$} is the set of all functions $f\colon\N\to\Fin{X}$.
This set is partially ordered by defining
\IndexPrint{$f\sub g$} as $f(n)\sub g(n)$ for all $n$.
\edfn

The cardinal $\cof(\Fin{X}^\N)$ depends only on $\card{X}$.

\blem\label{ukbfinkw}
Let $(G,T)$ be a metrizable boundedness system, and let $\kappa=\card{T}$. Then:
\be
\itm $\Fin{\kappa}^\N\preceq \Bdd_T(G)$.
\itm $\cof(\Bdd_T(G))\le\cof(\Fin{\kappa}^\N)$.
\itm $\cof(\PK(G))\le\cof(\Fin{\density(G)}^\N)$.
\ee
\elem
\bpf[Proof of (1)]
Fix a neighborhood base $U_n$, $n\in\N$, at $e$.
For each $f\in \Fin{\kappa}^\N$, define
$$K_f=\bigcap_{n\in\N} f(n)\mult U_n.$$
Then each set $K_f$ is in $\Bdd_T(G)$, and the family $\set{K_f}{f\in \Fin{\kappa}^\N}$ is cofinal in $\Bdd_T(G)$.
\epf

The following concept is central for the main results of this section.

\bdfn
The \emph{local density}\index{$\localdensity(G)$} of a group $G$ is the cardinal
$$\localdensity(G) = \min\set{\density(U)}{U\mbox{ is a neighborhood of }e\mbox{ in }G}.$$
$G$ has \emph{stable density} if $\localdensity(G)=\density(G)$.
\edfn

$G$ has local density $\kappa$ if, and only if,
$G$ has a local base at $e$, consisting of elements of density $\kappa$.

\blem
The cardinal
$\localdensity(G)$ is the minimal density of a clopen subgroup $H$ of $G$. Thus, $G$ has stable density if, and only if,
$\density(H)=\density(G)$ for all clopen $H\le G$.
\elem
\bpf
Let $U\sub G$ be an open neighborhood of $e$ with $\density(U)=\localdensity(G)$. Take $H=\la U\ra$.
Then $H$ is an open, and thus closed, subgroup of $G$.
\epf

\bexm
Every connected group has stable density.
\eexm

\bdfn
Let $V$ be a neighborhood of $e$ in $G$.
A set $A\sub G$ is a \emph{$V$-grid} if the sets $aV$, for $a\in A$, are pairwise disjoint.
A set $A$ is a \emph{grid} if it is a $V$-grid for some neighborhood $V$ of $e$.
\edfn

The intersection of a precompact set and a grid must be finite.

\blem\label{closeddiscrete}
Let $G$ be a metrizable group with stable density.
Let $\kappa=\density(G)$, and $U$ be a neighborhood of $e$.
\be
\itm For each $\lambda<\kappa$, the neighborhood $U$ contains a grid of cardinality
greater than $\lambda$.
\itm If $\cof(\kappa)>\alephes$, then $U$ contains a grid of cardinality $\kappa$.
\ee
\elem
\bpf
(1) Let $V\sub U$ be a symmetric neighborhood of $e$, such that for each $S\sub G$ with $\card{S}=\lambda<\kappa$,
$SV^2$ does not contain $U$.

By Zorn's Lemma, there is a maximal $V$-grid $A$ in $U$.
As $V$ is symmetric, $U\sub AV^2$.
It follows that $\card{A}>\lambda$.

(2) Let $\set{V_n}{n\in\N}$ be a symmetric local base at $e$, and for each $n$ let $A_n$ be a maximal $V_n$-grid in $U$.
The previous argument shows that for each $\lambda<\kappa$, there is $n$ such that $\card{A_n}>\lambda$.
Thus, $\sup_n\card{A_n}=\kappa$. As $\cof(\kappa)>\alephes$, there is $n$ with $\card{A_n}=\kappa$.
\epf

We are now ready for the main results of this section.
Given partially ordered sets $P_1,\dots,P_k$, define the \emph{coordinate-wise partial order} on $P_1\x \dots\x P_k$
by $(a_1,\dots,a_k)\le (b_1,\dots,b_k)$ if $a_1\le b_1$, \dots, $a_k\le b_k$.

\bdfn
For cardinals $\kappa,\lambda$, the family
$$[\kappa]^\lambda:=\set{A\sub\kappa}{\card{A}=\lambda}$$
is partially ordered by $\sub$.
\edfn

\bthm\label{goodthm}
Let $G$ be a metrizable non-locally precompact group of stable density $\kappa$.
Then $\cof(\PK(G))=\fd\cdot\cof([\kappa]^\alephes)$.
\ethm

Theorem \ref{goodthm} follows from the following two propositions.

\bprp\label{goodprp}
Let $G$ be a metrizable non-locally precompact group of stable density $\kappa$.
Then:
\be
\itm $\PK(G)$ is cofinally equivalent to $\Fin{\kappa}^\N$.
\itm $\cof(\PK(G))=\cof(\Fin{\kappa}^\N)$.
\ee
\eprp
\bpf
If $\cof(\kappa)>\alephes$, let $\kappa_n=\kappa$ for all $n$.
Otherwise, for $n\in\N$ let $\kappa_n$ be cardinals such that $\kappa_n<\kappa_{n+1}$ for all $n$
and $\sup_n\kappa_n=\kappa$.

Let $\set{U_n}{n\in\N}$ be a decreasing local base at $e$.
For each $n$, there is by Lemma \ref{closeddiscrete} a grid $A_n\sub U_n$
with $\card{A_n}=\kappa_n$.
Let $P\in\PK(G)$. Then $P\cap A_n$ is finite for all $n$.
Thus, we can define $\Psi\colon\PK(G)\to\prod_n\Fin{A_n}$ by
$$P\mapsto f\mbox{ with }f(n)=P\cap A_n$$
for all $n$.

The map $\Psi$ is cofinal:
For each $f\in \prod_n\Fin{A_n}$, the set $P=\Union_n f(n)\cup\{e\}$ is a countable
set converging to $e$, and thus compact. For each $n$, we have that $f(n)\sub \Psi(P)(n)$.

As $\Psi$ is monotone and cofinal, $\PK(G)\preceq\prod_n\Fin{A_n}$.

\blem
If $\kappa_n\le\kappa_{n+1}$ for all $n$, and $\sup_n\kappa_n=\kappa$, then
$$\prod_n\Fin{\kappa_n}\preceq\NN\x \prod_n\Fin{\kappa_n}\preceq\Fin{\kappa}^\N.$$
\elem
\bpf
To prove the first assertion, map $f$ to the pair $(h,f)$, where
$h(n)=\max f(n)\cap \w$ (or $0$ if $f(n)\cap \w$ is empty).
For the second assertion, map $(h,g)$ to the function
$f(n)=\Union_{m\le h(n)}g(m)$.
\epf

Finally, apply Lemma \ref{ukbfinkw}.
\epf

\bprp\label{morph}
For each infinite cardinal $\kappa$,
$\cof(\kfinw)=\fd\cdot\cof([\kappa]^\alephes)$.
\eprp
\bpf
$\kfinw\preceq\NN\x[\kappa]^\alephes$:
Given a function $f\in\kfinw$, define $g_f\in\NN$ by $g_f(n)=\max (f(n)\cap\N)\cup\{0\}$,
and $s_f=\Union_n f(n)$. The map $f\mapsto(g_f,s_f)$ is monotone and cofinal.
Thus, $\cof(\kfinw)\ge \cof(\NN\x[\kappa]^\alephes) = \fd\cdot\cof([\kappa]^\alephes)$.

$(\le)$ For each $s\in[\kappa]^\alephes$, fix a surjection $r_s\colon\N\to s$.
Consider the mapping of $(f,s)\in\NN\x[\kappa]^\alephes$ to $g\in\kfinw$,
defined by
$$g(n) = \{r_s(1),r_s(2),\dots,r_s(f(n))\}$$
for all $n$.
Then the image of a product of two cofinal sets is cofinal.
\epf

We now treat the general case, using the following observation:
If $H$ is a clopen subgroup of $G$ of density $\localdensity(G)$, then $H$ has stable density,
$G/H$ is discrete, and $\density(G)=\card{G/H}\cdot \localdensity(G)$.

\bthm\label{goodnews}
Let $G$ be a metrizable non-locally precompact group.
\be
\itm Let $H$ be a clopen subgroup of $G$, of density $\localdensity(G)$. Then
$\PK(G)$ is cofinally equivalent to $\Fin{G/H}\x\Fin{\localdensity(G)}^\N$.
\itm $\cof(\PK(G))=\fd\cdot \density(G)\cdot\cof([\localdensity(G)]^\alephes)$.
\ee
\ethm
\bpf
(1) $\density(H)=\localdensity(G)=\localdensity(H)$.

\blem
For each clopen subgroup $H$ of $G$, $\PK(G)$ is cofinally equivalent to $\Fin{G/H}\x \PK(H)$.
\elem
\bpf
Fix a set $S\sub G$ of coset representatives, that is such that $\card{S\cap gH}=1$ for all $g\in G$.
We need to show that $\PK(G)$ is cofinally equivalent to $\Fin{S}\x \PK(H)$.
For $A\sub G$ let $S(A)=\set{s\in S}{sH\cap A\neq\emptyset}$.
The function
$$P\mapsto \left(S(P),H\cap\Union_{s\in S(P)}s\inv P\right)$$
is a monotone and cofinal map from $\PK(G)$ to $\Fin{S}\x \PK(H)$.

For the other direction, we can map each $(F,P)\in \Fin{S}\x \PK(H)$ to $FP$.
\epf
This, together with Theorem \ref{goodthm}, proves (1).

(2) By (1) and Proposition \ref{morph},
$$\cof(\PK(G))=\card{G/H}\cdot\fd\cdot \cof([\localdensity(G)]^\alephes).$$
The statement follows, using that $\card{G/H}\le \density(G)\le \cof(\PK(G))$
(Lemma \ref{kappalecofbdd}).
\epf

\bexm
For all cardinals $\lambda\le\kappa$, there are metrizable groups $G$ with
$\localdensity(G)=\lambda$ and $\density(G)=\kappa$. For example, a product of a discrete group of cardinality
$\kappa$ and $C(\bbT^\lambda,\bbT)$.
An extreme example is where $G$ is discrete: We obtain $\localdensity(G)=1$, and $\density(G)=\card{G}$,
and indeed $\PK(G)=\Fin{G/\{e\}}$.
\eexm

The cardinal $\cof(\Fin{\kappa}^\N)$ also appears, in a different context, in a completely different context
studied by Bonanzinga and Matveev \cite{MilenaMisha}.

\section{Abelian groups and Pontryagin--van Kampen duality}\label{pvk}

In the remainder of the paper, all considered groups are assumed to be abelian,
and we use the additive notation and $0$ for the trivial element.
In particular, we identify $\bbT$ with the additive group $[-1/2,1/2)$,
having addition defined by identifying $\pm 1/2$.

\emph{A character}\index{character!a character on a topological abelian group} on a topological abelian group $G$ is a continuous group homomorphism from $G$ to
the torus group $\bbT$. This is a collision in terminology,
which may be solved as follows:
Characters on $G$ are its continuous homomorphisms into $\bbT$, whereas \emph{the} character of $G$ is the minimal
cardinality of a local base of $G$ at $e$.
The set of all characters on $G$, with pointwise addition, is a group.

For a topological abelian group $G$, let \IndexPrint{$\K(G)$} denote the family of all compact subsets of $G$.
For a set $A\sub G$ and a positive real $\epsilon$, define
$$[A,\epsilon] := \set{\chi \in \Gh{}}{\card{\chi(a)}\leq\epsilon\mbox{ for all }a \in A}.$$
The sets $[K,\epsilon]\sub\Gh{}$, for $K\in\K(G)$ and $\epsilon>0$, form a neighborhood base at the trivial character,
defining the compact-open topology.
We write \IndexPrint{$\Gh{}$} for the topological abelian group obtained in this manner.

A topological abelian group $G$ is \emph{reflexive} if the evaluation map\index{evaluation map, $E$}
$$E\colon G\to\widehat{\Gh{}},$$
defined by $E(g)(\chi)=\chi(g)$ for all $g\in G$ and $\chi\in\Gh{}$,
is a topological isomorphism.
By the Pontryagin--van Kampen theory, we know that
every locally compact abelian group is reflexive.
Furthermore, the dual of a compact group is discrete and the dual
of a discrete group is compact.
In general, the dual of a locally compact abelian group is also locally compact.
It follows that every compact abelian group is equipped with
the topology of pointwise convergence on its dual group.
This fact will be used below.

Let $K$ be a compact subset of $G$. For each $n$, the set $K_n=K\cup 2K\cup\dots\cup nK$ is
compact, and $[K_n,1/4]\sub [K,1/4n]$. Thus, the sets $[K,1/4]$, for $K\in\K(G)$, also form a neighborhood
base of $\widehat G$ at the trivial character.

\bdfn
Let $G$ be a topological abelian group.
For $A\sub G$, let $A^\rhd := [A,1/4]$\index{$A^\rhd$}.
Similarly, for $X\subseteq\Gh{}$, let
$$X^\lhd :=\Bigl\{\, g\in G : \card{\chi(g)}\leq \frac{1}{4}\mbox{ for all }\chi \in X\,\Bigr\}.$$
\edfn

\blem[{\cite[Proposition 1.5]{banasz:nuclear}}]\label{Upol}
For each neighborhood $U$ of $0$ in $G$, we have that $U^\rhd\in\K(\Gh)$.
\elem

\bdfn[Vilenkin \cite{vilenkin}]
Let $G$ be a topological abelian group.
A set $A\sub G$ is \emph{quasiconvex} if $A^{\rhd\lhd}=A$.
The topological group
$G$ is \emph{locally quasiconvex} if it has a neighborhood base at its identity, consisting of quasiconvex sets.
\edfn

For each set $A\sub G$, the set $A^\rhd$ is a quasiconvex subset of $\Gh{}$.
Thus, the topological group $\Gh{}$ is locally quasiconvex for all topological abelian groups $G$.
Moreover, local quasiconvexity is hereditary for arbitrary subgroups.

The set $A^{\rhd\lhd}$ is the smallest quasiconvex subset of $G$ containing $A$. This set is closed.

In the case where $G$ is a topological vector space $G$ is locally quasiconvex in the present sense
if, and only if, $G$ is a locally convex topological vector space in the ordinary sense \cite{banasz:nuclear}.

If $G$ is locally quasiconvex, its characters separate points of $G$,
and thus the evaluation map $E\colon G\to \Ghh{}$ is injective.
For each quasiconvex neighborhood $U$ of $0$ in $G$, the set $U^\rhd$ is a compact subset of $\Gh{}$
(Lemma \ref{Upol}), and thus $U^{\rhd\rhd}$ is a neighborhood of $0$ in $\Ghh{}$.
As $E[G]\cap U^{\rhd\rhd}=E[U^{\rhd\lhd}]=E[U]$, we have that $E$ is open \cite[Lemma 14.3]{banasz:nuclear}.

\blem\label{cofeq1}
Let $G$ be a complete locally quasiconvex group.
Let \IndexPrint{$\Nh$} be the family of all neighborhoods of $0$ in $\Gh$.
Then:
\be
\itm $(\Nh,\spst)$ is cofinally equivalent to $(\K(G),\sub)$.
\itm $\chi(\Gh)=\cof(\K(G))$.
\ee
\elem
\bpf[Proof of (1)]
We have seen above that the monotone map $\rhd\colon \K(G)\to\Nh$ is cofinal.

Consider the other direction. Let $K\in\K(G)$, and take $U=K^\rhd\in\Nh$.
By Lemma \ref{Upol}, $U^\rhd\in\K(\Ghh)$. Now,
$$K\sub K^{\rhd\lhd}=U^\lhd=E\inv[U^\rhd\cap E[G]].$$
As $G$ is complete, $U^\rhd\cap E[G]$ is compact. As $G$ is locally quasiconvex, the map $E$ is open, and therefore $E\inv[U^\rhd\cap E[G]]$
is compact. Thus, the monotone map $\lhd\colon\Nh\to\K(G)$ is also cofinal.
\epf

\brem
As can be seen from the proof of Lemma \ref{cofeq1}, the assumption that $G$ is complete can be wakened to the so-called
\emph{quasiconvex compactness property}, that is, the property that for each $K\in\K(G)$, we have that
$K^{\rhd\lhd}\in\K(G)$.
\erem

We obtain the following proposition, which extends to topological abelian groups a
result of Saxon and Sanchez--Ruiz for the strong dual of a
metrizable space \cite[Corollary 2]{SaxSan95}.
As every locally convex topological
vector space is connected, it has stable density and therefore the concept of local
density is not required in \cite{SaxSan95}.
As stated here, our result does not generalize that of Saxon and Sanchez--Ruiz. There
is a natural extension of our approach which implies their result as well,
by replacing $\K(G)$ with more general boundedness notions on $G$. For concreteness, we do not
present our results in full generality.

A topological space $X$ is a \emph{$k$-space} if the
topology of $X$ is determined by its compact subsets, that is,
$F\sub X$ is closed if (and only if) $F\cap K$ is
closed in $K$ for all $K\in\K(G)$. Every metrizable space is a $k$-space.
A \emph{$k$-group} is a topological group which is a $k$-space.

Let $G$ be the dual of a metrizable group $\Gamma$. If $\Gamma$ is (pre)compact,
then by Pontryagin's Theorem, $G$ is discrete, that is $\chi(G)=1$.
Item (1) of the following proposition is known \cite[Theorem 3.12(ii)]{ComfortHBK}.

\bprp\label{prop_d_of_M}
Let $G$ be the dual of a metrizable, nonprecompact group $\Gamma$.
\be
\itm If $\Gamma$ is locally precompact, then $\chi(G)=\density(\Gamma)$.
\itm If $\Gamma$ is non-locally precompact,
then $\chi(G)$ is the maximum of $\fd$, $\density(\Gamma)$, and $\cof([\localdensity(\Gamma)]^\alephes)$.
\ee
\eprp
\bpf
Au\ss{}enhofer \cite{aus} and, independently, Chasco \cite{Chas:97} proved that
a metrizable group and its completion have the same (topological) dual group.
Since the density and local density of a metrizable group are equal to those
of its completion, we may assume that $\Gamma$ is complete.

Since $\Gamma$ is metrizable, it is a $k$-space, and therefore $G=\widehat{\Gamma}$ is complete
\cite[Proposition 1.11]{banasz:nuclear}.
By Lemma \ref{cofeq1} and the completeness of $\Gamma$, we have that
$$\chi(G)=\chi(\widehat\Gamma)=\cof(\K(\Gamma))=\cof(\PK(\Gamma)).$$

(1) By Lemma \ref{lb}, $\cof(\PK(\Gamma))=\density(\Gamma)$.

(2) By Theorem \ref{goodnews} and Theorem \ref{morph}, we have that
$$\cof(\PK(\Gamma))=\density(\Gamma)\cdot\cof(\Fin{\localdensity(\Gamma)}^\N)=
\fd\cdot \density(\Gamma)\cdot\cof([\localdensity(\Gamma)]^\alephes).\qedhere$$
\epf

Even for locally quasiconvex $G$, the evaluation map $E$ need not be continuous.
If it is, then $G$ is isomorphic to its image $E[G]$ in $\Ghh{}$.

\bdfn
A topological abelian group $G$ is \emph{subreflexive} if the evaluation map $E\colon G\to E[G]$ is a topological isomorphism.
In this case, we identify $G$ with its image $E[G]\le \Ghh{}$.
\edfn

\brem
If $G$ is a subreflexive topological abelian group, then $G$ is locally quasiconvex.
Indeed, the group $\Ghh$ is locally quasiconvex, being a dual group, and
therefore so is its subgroup $E[G]$, which is isomorphic to $G$.
\erem

\blem\label{polar2}
Let $G$ be a subreflexive topological abelian group.
Then the family $\set{K^\lhd}{K\in\K(\Gh{})}$ is a neighborhood base at $e$ in $G$.
\elem
\bpf
Let $K\in\K(\Gh{})$. The set $K^\rhd$ is a neighborhood of $0$ in $\Ghh{}$.
As $G$ is subreflexive, $K^{\lhd}$ is a neighborhood of $0$ in $G$.

Let $U$ be a neighborhood of $e$ in $G$. As $G$ is locally quasiconvex, we may assume that $U$ is quasiconvex.
Then the set $K:=U^{\rhd}$ is compact in $\Gh$ (Lemma \ref{Upol}), and $K^\lhd=U^{\rhd\lhd}=U$.
\epf


\bprp\label{cofeq}
Let $G$ be a subreflexive topological abelian group, and $\cN$ be the family of all neighborhoods of $0$ in $G$.
Then:
\be
\itm $(\cN,\spst)$ is cofinally equivalent to $(\K(\Gh{}),\sub)$.
\itm $\chi(G)=\cof(\K(\Gh{}))$.
\ee
\eprp
\bpf[Proof of (1)]
By Lemma \ref{polar2}, the monotone map $\lhd\colon\K(\Gh{})\to\cN$ is cofinal.
The monotone map $\rhd\colon\cN\to\K(\Gh{})$ is also cofinal: Let $K\in\K(\Gh{})$.
By Lemma \ref{polar2}, $K^\lhd\in\cN$, and $(K^\lhd)^\rhd\spst K$.
\epf

Even complete subreflexive groups $G$ need not be reflexive.
The following corollary tells that, however, $\Ghh{}$ is not much larger than $G$.
(See also Theorem \ref{bb} and Corollary \ref{ww} below.)
Au\ss{}enhofer made related observations in \cite[5.22]{aus}.
Question 5.23 in \cite{aus} asks whether the character group of an abelian
metrizable group is reflexive.

\bcor\label{chighh}
\mbox{}
\be
\itm For subreflexive $G$ with $\Gh$ complete, $\chi(\Ghh)=\chi(G)$.
\itm If $G$ is a locally quasiconvex $k$-group, then $\chi(\Ghh)=\chi(G)$.
\ee
\ecor
\bpf
(1) The group $\Gh$ is locally quasiconvex. By Lemma \ref{cofeq1} and Proposition \ref{cofeq}, $\chi(\Ghh)=\cof(\K(\Gh))=\chi(G)$.

(2) By Corollary \ref{lcqpkiswr} below, the group $G$ is subreflexive.
As $G$ is a $k$-group, the group $\Gh$ is complete. Apply (1).
\epf

The first two items in the following theorem are well known.

\bthm\label{Ghmetric}
Let $G$ be a subreflexive group such that the group $\Gamma=\Gh$ is metrizable.
Then $\chi(G)=\cof(\PK(\Gamma))$. Thus,
\be
\itm If $\Gamma$ is precompact, then $\chi(G)=1$, that is, the topological group $G$ is discrete.
\itm If $\Gamma$ is nonprecompact locally precompact, then $\chi(G)=\density(\Gamma)$.
\itm If $\Gamma$ is non-locally precompact, then $\chi(G)=\fd\cdot \density(\Gamma)\cdot\cof([\localdensity(\Gamma)]^\alephes)$.
\ee
\ethm
\bpf
By Proposition \ref{cofeq}, we have that $\chi(G)=\cof(\K(\Gh{}))=\cof(\K(\Gamma))$.
Let $\Delta$ be the completion of $\Gamma$. The group $\Delta$ is locally quasiconvex too, and metrizable,
and thus subreflexive.
By Corollary \ref{completion}, we have that $\cof(\K(\Delta))=\cof(\PK(\Gamma))$.

It remains to prove that $\K(\Gamma)$ is cofinally equivalent to $\K(\Delta)$.
By the Au\ss{}enhofer--Chasco Theorem, we may identify $\widehat\Delta$ with $\widehat\Gamma$.
As $G$ is subreflexive, we may also identify $G$ with its image in $\Ghh=\widehat\Gamma$,
and similarly for $\Delta$.

$\K(\Delta)\preceq \K(\Gamma)$: Let $K\in\K(\Delta)$.
Then $K^\rhd$ is a neighborhood of $0$ in $\widehat\Delta=\widehat\Gamma=\Ghh$.
As $G$ is subreflexive, $K^\rhd\cap G$ is a neighborhood of $0$ in $G$, and thus $(K^\rhd\cap G)^\rhd\in\K(\Gh)=\K(\Gamma)$.
Define $\Phi(K)=(K^\rhd\cap G)^\rhd$.
For each $K\in\K(\Gamma)$, $K\in\K(\Delta)$ and $\Phi(K)\spst K$. Thus, $\Phi$ is cofinal.

$\K(\Gamma)\preceq\K(\Delta)$:
Let $K\in\K(\Gamma)$. Then $K^\rhd$ is a neighborhood of $0$ in $\widehat\Gamma=\widehat\Delta$.
Thus, $K^{\rhd\rhd}\in\K(\Delta\hat{~}\hat{~})$, and as $\Delta$ is complete,
$K^{\rhd\rhd}\cap\Delta\in\K(\Delta)$.
Define $\Psi\colon\K(\Gamma)\to\K(\Delta)$ by $\Psi(K)=K^{\rhd\rhd}\cap\Delta$.
For each $C\in \K(\Delta)$, $C^\rhd$ is a neighborhood of $0$ in $\widehat\Delta=\widehat\Gamma$,
and thus there is $K\in\K(\Gamma)$ such that $K^\rhd\sub C^\rhd$.
Then $K^{\rhd\rhd}\spst C^{\rhd\rhd}\spst C$, and therefore $\Psi(K)=K^{\rhd\rhd}\cap\Delta\spst C$.
This shows that $\Psi$ is cofinal.

(1) and (2) follow, using Lemma \ref{lb} and Theorem \ref{goodnews}.
\epf

Theorem \ref{Ghmetric} is stronger than Proposition \ref{prop_d_of_M}: duals of metrizable groups
are subreflexive, and have a metrizable dual.

\section{Application to the free abelian topological groups}\label{AXsec}

A topological space $X$ is \emph{hemicompact} if $\cof(\K(X))\le\alephes$.
$X$ is a \emph{$k_\w$ space} if it is a hemicompact $k$-space.
Denote the weight of a topological space $X$ by \IndexPrint{$\weight(X)$}.

The following theorem extends, but does not generalize,
several results of Nickolas and Tkachenko (e.g., the results numbered 2.12, 2.18, 2.22 in \cite{NT1},
and those numbered 2.9, 3.5, 3.7 in \cite{NT2}.)
For example, Nickolas and Tkachenko proved that if $X$ is \emph{compact}, then
$$\chi(A(X))=\fd\cdot\cof([\weight(X)]^\alephes),$$
and that if $X$ is a $k_\w$ space such that all compact subsets of $X$ are metrizable, then
$\chi(A(X))=\fd$.
Nickolas and Tkachenko's results were proved by direct, but more combinatorially involved,
methods.

\bthm\label{chiAX}
Let $X$ be a nondiscrete $k_\w$ space of compact weight $\kappa$. Then
$$\chi(A(X))=\fd\cdot\cof([\kappa]^\alephes).$$
\ethm
\bpf
Au\ss{}enhofer \cite{aus} and, independently, Galindo--Hern\'andez \cite{GalindoHernandez99} proved
that for a class of spaces $X$ containing $k$-spaces (namely, Ascoli $\mu$-spaces),
the free abelian topological group $A(X)$ is subreflexive.
Pestov \cite{Pestov95} proved that for a class of spaces $X$ containing $k_\w$ spaces (namely, $\mu$-spaces),
$\widehat{A(X)}= C(X,\bbT)$.
As $X$ is $k_\w$, $C(X,\bbT)$ has a countable local base at $0$
(namely, the sets $[K_n,1/n]$ where $\set{K_n}{n\in\N}$ is cofinal in $\K(X)$).
Thus, $C(X,\bbT)$ is metrizable.

Moreover, $C(X,\bbT)$ is
non-locally precompact. Thus, Theorem \ref{Ghmetric} applies.

\blem\label{CXTkw}
Let $X$ be a Tychonoff space of compact weight $\kappa$. Then:
\be
\itm $\boundedness(C(X,\bbT))=\boundedness(C(X,\R))=\kappa$.
\itm If $X$ is hemicompact (or just $\cof(K(X))\le\kappa$), then
$$\boundedness(C(X,\bbT))=\density(C(X,\bbT))=\localdensity(C(X,\bbT))= \weight(C(X,\bbT))=\kappa.$$
In particular, $C(X,\bbT)$ has stable density.
\ee
\elem
\bpf
For each cofinal family $\cK\sub \K(X)$, and for $Y=\bbT$ or $\R$,
the mapping $f\mapsto (\, f|_K : K\in\cK \,)$ is an embedding of $C(X,Y)$
in $\prod_{K\in\cK}C(K,Y)$.

(1)
If $X$ is locally compact and $\weight(X)$ is infinite, then
$\weight(C(X,\bbT))\le \weight(X)$
\cite[3.4.16]{Engelking}.
Thus, in the case $\cK=\K(X)$, we have that
\begin{eqnarray*}
\boundedness(C(X,Y)) & \le & \boundedness\bigl(\prod_{K\in\K(X)}C(K,Y)\bigr)=\sup_{K\in\K(X)} \weight(C(K,Y))\le\\
& \le &  \sup_{K\in\K(X)} \weight(K),
\end{eqnarray*}
Let $K\in\K(X)$. Take $S\sub C(X,Y)$ with $\card{S}=\boundedness(C(X,Y))$, such that $S+[K,1/16]=C(X,Y)$.
Then $\set{f\inv(-1/16,1/16)\cap K}{f\in S}$ is a base of $K$:
Let $p\in U\cap K$, $U$ open in $X$. As $X$ is Tychonoff, there is $g\in C(X,Y)$ such
that $g$ is $1/4$ on $X\sm U$ and $g(p)=0$. As $S+[K,1/16]=C(X,Y)$, there is
$f\in S$ such that $\card{f(x)-g(x)}\le 1/16$ for each $x\in K$. It follows that
$p\in g\inv(-1/16,1/16)\cap K\sub U\cap K$.
Thus, $\weight(K)\le \boundedness(C(X,Y))$ for each $K\in\K(X)$.

(2) By (1), $\kappa=\boundedness(C(X,\R))\le \density(C(X,\R))$. As $C(X,\R)$ is connected, $\density(C(X,\R))=\localdensity(C(X,\R))$.
For each $\epsilon<1/2$ and each compact $K\sub X$, $[K,\epsilon]$ is the same in $C(X,\R)$
and in $C(X,\bbT)$.
Thus,
$$\kappa\le \localdensity(C(X,\R))\le \localdensity(C(X,\bbT))\le \density(C(X,\bbT))\le \weight(C(X,\bbT)).$$
In the case where $\card{\cK}=\cof(K(X))$,
\begin{eqnarray*}
\weight(C(X,\bbT)) & \le & \weight\left(\prod_{K\in\cK}C(K,\bbT)\right)=\card{\cK}\cdot\sup_{K\in\cK} \weight(C(K,\bbT)) \le\\
& \le & \cof(\K(X))\cdot \sup_{K\in\K(X)} \weight(K)\le\kappa\cdot\kappa=\kappa.\qedhere
\end{eqnarray*}
\epf

We therefore have, by Theorem \ref{Ghmetric}, that $\chi(A(X))$ is the maximum of
$\fd$ and $\cof([\kappa]^\alephes)$, where $\kappa=\density(C(X,\bbT))=\sup\set{\weight(K)}{K\in\K(X)}$.
This completes the proof of Theorem \ref{chiAX}.
\epf

\bexm
If $X$ is compact, or locally compact $\sigma$-compact,
then $X$ is a $k_\w$ space and Theorem \ref{chiAX} applies.
\eexm

As already pointed out in the introduction, by virtue of \cite[Corollary 2.3]{NT2}, 
our results also apply to the free \emph{nonabelian} topological groups $F(X)$.

\section{The inner theorem}

We begin with an inner characterization of subreflexivity.

\bdfn
A set $V\sub G$ is a \emph{$k$-neighborhood} of $0$ if
for each $K\in\K(G)$ with $0\in K$, $V\cap K$ is a neighborhood of $0$ in $K$.
\edfn

\blem[Hern\'{a}ndez--Trigos--Arrieta \cite{her_tri:05}]\label{le_pk_group}
\mbox{}
\be
\itm Let $G$ be a $k$-group.
Every quasiconvex $k$-neighborhood of $0$ is a neighborhood of $0$.
\itm Let $U$ be a quasiconvex subset of a locally quasiconvex group
$G$. $U$ is a $k$-neighborhood of $0$ if, and only if, $U^\rhd\in\K(\Gh{})$.
\ee
\elem

We obtain the following.

\bthm\label{th_polar}
A group $G$ is subreflexive if, and only
if, $G$ is locally quasiconvex, and each quasiconvex $k$-neighborhood of the identity in $G$ is a neighborhood of the identity.
\ethm
\bpf
$(\Leftarrow)$ Let $F\in\K(\widehat G)$ and $K\in\K(G)$.
By Ascoli's Theorem, the restrictions of the elements of $F$ to $K$ form
an equicontinuous subset of $C(K,\bbT)$. Hence, if $K$ contains
$0$, then $F^\rhd\cap K$
is a neighborhood of $0$ in $K$.
Again, taking intersections, we have that $F^\lhd\cap K$
is a neighborhood of $0$ in $K$. Thus, $F^\lhd$ is a
neighborhood of $0$.

$(\Impl)$ Let $W$ be a quasiconvex $k$-neighborhood of $0$. Then $W^\rhd$ is compact in $\Gh{}$.
As $G$ is subreflexive, $W=W^{\rhd\lhd}$ is a neighborhood of $0$ in $G$.
\epf

Lemma \ref{le_pk_group} and Theorem \ref{th_polar} imply the following.

\bcor[folklore]\label{lcqpkiswr}
Every locally quasiconvex $k$-group is subreflexive.\qed
\ecor

For locally convex topological vector spaces and countable weight,
the following result was proved by Ferrando, K\c{}akol, and M. L\'opez Pellicer
\cite{fer_kak_lop:06}.

\bthm\label{dwt}
Let $G$ be a locally quasiconvex abelian group.
\be
\itm The cardinal $\boundedness(\Gh{})$ is equal to the compact weight of $G$.
\itm If the topological group $\Gh{}$ is metrizable, then $\density(\Gh{})$ equal to the compact weight of $G$.
\ee
\ethm
\bpf[Proof of (1)]
$(\le)$
As $\Gh\le C(G,\bbT)$, we have by Lemmata \ref{subgp-kbdd} and \ref{CXTkw}
that $\boundedness(\Gh)\le \boundedness(C(G,\bbT))=\sup\set{\weight(K)}{K\in\K(G)}$.

$(\ge)$ Let $K\in\K(G)$. Since $[K,1/8]$ is a neighborhood of the identity of $\Gh{}$, there is a set $S\sub\Gh{}$
with $\card{S}\le \boundedness(\Gh{})$ such that $S+[K,1/8]=\Gh{}$.

The set $S$ separates the points of $K$:
Let $a_1,a_2$ be distinct elements of $K$.
As $G$ is locally quasiconvex, there is $\chi \in \Gh$ such that
$\card{\chi(a_1-a_2)}>1/4$. As $\chi\in\Gh=S+[K,1/8]$,
there are $\alpha\in S$ and $\beta\in [K,1/8]$ such that
$\chi=\alpha+\beta$. Then $\card{\beta(a_1-a_2)}\leq \card{\beta(a_1)}+\card{\beta(a_2)}\le 2/8=1/4$, and
thus $\card{\alpha(a_1-a_2)}\geq \card{\chi(a_1-a_2)}-1/4>0$.

Thus, the minimal topology on $K$ which makes all elements of $S$ continuous is Hausdorff,
and as $K$ is compact, its topology (which is minimal Hausdorff) coincides with it.
Thus, $\weight(K)\leq \card{S}\le \boundedness(\Gh{})$.
\epf

An unpublished result of Au\ss{}enhofer asserts that,
if $G$ is a separable metrizable group, then all higher character groups of $G$
are separable. This is in accordance with item (3) of the following theorem.

\bthm\label{bb}
Let $G$ be a topological abelian group, and let $\kappa$ be the compact weight of $\Gh$.
\be
\itm If $G$ is subreflexive then $\boundedness(G)=\boundedness(\Ghh{})=\kappa$.
\itm If $G$ is a locally quasiconvex $k$-group then $\boundedness(G)=\boundedness(\Ghh{})=\kappa$.
\itm If $G$ is locally quasiconvex and metrizable then $\density(G)=\density(\Ghh{})=\kappa$.
\ee
\ethm
\bpf
(1) As $G\le\Ghh$, we have that $\boundedness(G)\le \boundedness(\Ghh)$.
By Theorem \ref{dwt}, $\boundedness(\Ghh)=\kappa$.
We prove that $\kappa\le \boundedness(G)$.

Let $K$ be a compact subset of $\Gh{}$. As $G$ is subreflexive, the set
$$U=(K\cup 2K)^\lhd=\set{g\in G}{(\forall \chi\in K)\ \card{\chi(g)}\le 1/8}$$
is a neighborhood of $0$ in $G$.
Let $S\sub G$ be such that $\card{S}\le \boundedness(G)$, and $S+U=G$.

$S$ separates points of $K$: Let $\chi,\psi\in K$ be distinct.
As $G^\rhd=\{0\}$, there is $g\in G$ such that $\card{(\chi-\psi)(g)}>1/4$.
Take $s\in S, u\in U$, such that $g=s+u$.
Then
$$\card{(\chi-\psi)(s)}\ge\card{(\chi-\psi)(g)}-\card{(\chi-\psi)(u)}>1/8.$$
It follows that $\weight(K)\le\card{S}\le \boundedness(G)$.

(2) Locally quasiconvex metrizable groups are subreflexive, being locally quasiconvex $k$-groups (Corollary \ref{lcqpkiswr}).
\epf

Mikhail Tkachenko pointed out to us that our results imply the following.

\bcor\label{ww}
For all subreflexive $G$ with $\Gh$ complete, $\weight(\Ghh)=\weight(G)$.
\ecor
\bpf
This follows from Corollary \ref{chighh} and Theorem \ref{bb}, using the
fact $\weight(G)=\boundedness(G)\cdot\chi(G)$ for all topological groups \cite{ArTkBk}.
\epf

We now turn to characterizing the local density of $\Gh$ in terms of inner properties of $G$.

A mapping is \emph{compact covering} if each compact subset of the range space is covered by
the image of a compact subset of the domain.

\blem\label{cc}
Let $H$ be a compact subgroup of $G$.
Then the canonical projection $\pi\colon G\to G/H$ is compact covering.
\elem
\bpf
For each compact $K\sub G/H$, the set $\pi\inv[K]$ is compact.
\epf

\blem\label{GoverK}
Let $G$ be a topological abelian group. Then:
\be
\itm For each compact subgroup $H$ of $G$, the topological groups $\widehat{G/H}$ and
$H^\rhd$ are isomorphic.

\itm For each open subgroup $H$ of $G$, the topological groups $\widehat{G/H}$
and $H^\rhd$ are isomorphic.
\ee
\elem
\bpf
(1) The homeomorphism $\vphi\colon\widehat{G/H}\to \Gh$ defined by $\vphi(\chi)=\chi\circ\pi$
is continuous and injective, and its image is $\set{\chi\in\Gh{}}{\chi|_H=0}=H^\rhd$.
A mapping is \emph{compact covering} if each compact subset of the range
space is covered by the image of a compact subset of the domain.
If $H$ be a compact subgroup of $G$,
then the canonical projection $\pi\colon G\to G/H$ is compact covering.

To see that $\vphi$ is open, let $U$ be a neighborhood of the identity of
$\widehat{G/H}$.
We may assume that $U=K^\rhd$ for some compact set $K\sub G/H$.
Since $\pi$ is compact covering, we may
assume that $K=\pi[K']$ for some compact set $K'\in\K(G)$.
We may also assume that $K'\spst H$.
Then $K'^\rhd\sub H^\rhd$, and therefore the set
$$\vphi[U]=\vphi[\pi[K']^\rhd]=\set{\vphi(\chi)}{\chi\in\pi[K']^\rhd}
=\set{\chi\circ\pi}{\chi\circ\pi\in K'^\rhd}=
K'^\rhd$$
is open.

(2) By the Pontryagin--van Kampen theorem,
since the group $G/H$ is discrete, the compact group $\widehat{G/H}$ is equipped
with the pointwise convergence topology.
As a consequence, the homeomorphism $\vphi\colon\widehat{G/H}\to \Gh$ defined by
$\vphi(\chi)=\chi\circ\pi$
is continuous and injective, and its image is $\set{\chi\in\Gh{}}{\chi|_H=0}=H^\rhd$.
The map $\vphi$ is also open since $\widehat{G/H}$ is compact.
\epf

For brevity, denote the compact weight of a group $G$ by \IndexPrint{$\kw(G)$}.

\bprp\label{hihi}
Let $G$ be a {locally quasiconvex} $k_\w$ group.
Then
$$\localdensity(\Gh)=\min\set{\kw(G/H)}{H\le G\mbox{ compact}}.$$
\eprp
\bpf
$(\ge)$ Let $\Gamma$ be an open subgroup of $G$ such that $\density(\Gamma)=\localdensity(\Gh)$.
As $G$ is $k_\w$, $\Gh$ is first countable and thus metrizable.
By Corollary \ref{lcqpkiswr}, the group $G$ is subreflexive.
As $k_\w$ groups are complete, $\Gamma^\lhd = \Gamma^\rhd\cap G$ is an intersection of a compact
group and a complete group, and is thus compact.

By Lemma \ref{GoverK}, $\widehat{G/\Gamma^\lhd}$ is isomorphic to $\Gamma^{\lhd\rhd}$, which contains $\Gamma$.
By definition, $\Gamma$ separates the points of $G/\Gamma^\lhd$, and therefore so does every dense subset
of $\Gamma$. Thus, $\weight(K)\le \density(\Gamma)$ for all compact sets $K\sub G/\Gamma^\lhd$.

$(\le)$ Let $H$ be a compact subgroup of $G$. By Lemma \ref{GoverK},
$\widehat{G/H}$ is isomorphic to $H^\rhd$. As $H^\rhd\le\Gh$, it is metrizable,
and thus by Corollary \ref{dwt},
$$\density(H^\rhd)=\density(\widehat{G/H})=\kw(G/H).$$
As $H^\rhd$ is open, $\localdensity(\Gh)\le \density(H^\rhd)$.
\epf

$G$ is \emph{locally hemicompact} (respectively, \emph{locally $k_\omega$})
if $G$ contains an \emph{open} hemicompact (respectively, $k_\omega$) subgroup.
The first item of the following theorem is an immediate consequence of the
Pontryagin--van Kampen theorem. The second item is new.

\bthm\label{th_hemi}
Let $G$ be a {locally quasiconvex}, locally $k_\w$ group. Let $H$ be an open $k_\w$ subgroup of $G$, of compact weight $\kappa$.
Let $\lambda = \min\set{\kw(H/K)}{K\le H\mbox{ compact}}$. Then:
\be
\itm If $H$ is nondiscrete and locally compact then $\chi(G)=\kappa$.
\itm If $H$ is non-locally compact then $\chi(G)$ is the maximum of $\fd$, $\kappa$ and $\cof([\lambda]^\alephes)$.
\ee
\ethm
\bpf[Proof of (2)]
As $H$ is open in $G$, $\chi(G)=\chi(H)$.
$G$ is locally quasiconvex, and therefore so is $H$.
By Lemma \ref{lcqpkiswr}, $H$ is subreflexive.
By hemicompactness, $\Gamma:=\widehat{H}$ is metrizable.
By Theorem \ref{Ghmetric},
$$\chi(H) = \fd\cdot \density(\Gamma)\cdot \cof([\localdensity(\Gamma)]^\alephes).$$
By Theorem \ref{dwt}(2), we have that $\density(\Gamma)=\kappa$.
By Proposition \ref{hihi}, $\localdensity(\Gamma)=\lambda$.
\epf

Concrete estimations are given in the overview (Section \ref{sec:overview}).
The proofs for these estimations are provided in the following, last section.

\section{Shelah's theory of possible cofinalities}\label{sec:pcf1}

In this section, we provide estimations for the cardinal $\cof([\kappa]^\alephes)$.
The estimations given here either appear explicitly in works of Shelah,
or are easy consequences thereof. Since we could not find a convenient reference for
these, we also provide proofs.

\blem\label{gtk}
For each cardinal $\kappa>\alephes$, we have that $\kappa\le\cof([\kappa]^\alephes)\le\kappa^\alephes$.
\elem
\bpf
Clearly, $\cof([\kappa]^\alephes)\le \card{[\kappa]^\alephes}=\kappa^\alephes$.
For the other inequality, note that if $A\sub [\kappa]^\alephes$ and $\card{A}<\kappa$,
then $\card{\Union A}\le \card{A}\cdot\alephes <\kappa$, and
thus $\Union A\neq\kappa$. In particular, $A$ is not cofinal in $[\kappa]^\alephes$.
\epf

For each cardinal $\lambda$,
the cardinal $\kappa=\lambda^\alephes$ has the property
$\kappa^\alephes=\kappa$.
This property holds for every cardinal $\kappa=2^\lambda$ for an infinite cardinal $\lambda$,
and if $\kappa^\alephes=\kappa$, then
the same is true for the subsequent cardinal $\kappa^+$.
This is also the case when $\kappa$ is inaccessible.
If the Generalized Continuum Hypothesis (GCH) holds, then this is the case for all cardinals
of uncountable cofinality.

\bcor\label{wpower}
For each infinite cardinal $\kappa$ with $\kappa^\alephes=\kappa$, we have that
$\cof(\kfinw)=\cof([\kappa]^\alephes)=\kappa$.
\ecor
\bpf
If $\kappa^\alephes=\kappa$, then $\kappa\ge \fc\ge\fd$. Apply Theorem \ref{morph} and Lemma \ref{gtk}.
\epf

\blem\label{bigcof}
For each $\kappa>\alephes$,
$\cof([\kappa]^\alephes)=\kappa\cdot\sup\set{\cof([\lambda]^\alephes)}{\lambda\le\kappa, \cof(\lambda)=\alephes}$.
\elem
\bpf
$(\ge)$ Monotonicity and Lemma \ref{gtk}.

$(\le)$ If $\cof(\kappa)=\alephes$, this follows from the fact that $\kappa\le\cof([\kappa]^\alephes)$ (Lemma \ref{gtk}).

If $\cof(\kappa)>\alephes$, then each countable subset of $\kappa$ is bounded in $\kappa$.
Thus, $[\kappa]^\alephes=\Union_{\alpha<\kappa}{[\alpha]^\alephes}$,
and therefore $\cof([\kappa]^\alephes)\le \kappa\cdot\sup\set{\cof([\lambda]^\alephes)}{\lambda<\kappa}$.
The statement for $\kappa=\aleph_1$ follows, and by induction, for each $\lambda<\kappa$ with $\lambda>\aleph_1$,
\begin{eqnarray*}
\cof([\lambda]^\alephes) & = & \lambda\cdot\sup\set{\cof([\mu]^\alephes)}{\mu\le\lambda, \cof(\mu)=\alephes}\le\\
& \le & \kappa\cdot\sup\set{\cof([\mu]^\alephes)}{\mu\le\kappa, \cof(\mu)=\alephes}.\qedhere
\end{eqnarray*}
\epf

\bcor\label{nojump0}
For each $\kappa$, if $\cof([\kappa]^\alephes)=\kappa$, then $\cof([\kappa^+]^\alephes)=\kappa^+$.\qed
\ecor

Item (1) of the following corollary is well known \cite{AbrMagHBK}, and Item (2) was proved, independently,
by Bonanzinga and Matveev \cite{MilenaMisha}.

\bcor\label{nojump}
\mbox{}
\be
\itm $\cof([\alephes]^\alephes)=1$, and for each $n\ge 1$, $\cof([\aleph_n]^\alephes)=\aleph_n$.
\itm $\cof(\finw{\alephes})=\fd$, and for each $n\ge 1$, $\cof(\finw{\aleph_n})=\fd\cdot\aleph_n$.\qed
\ee
\ecor

Already for $\kappa=\aleph_\omega$, the situation is different.
A diagonalization argument as in K\"onig's Lemma shows that,
$\cof([\kappa]^{\cof(\kappa)})>\kappa$ for singular cardinals $\kappa$.

\bcor
If $\cof(\kappa)=\alephes<\kappa$, then $\cof(\kfinw)\ge\fd\cdot\kappa^+$.\qed
\ecor

We next consider upper bounds.

\subsection{In the absence of large cardinals}

\emph{Shelah's Strong Hypothesis (SSH)} is the statement that
for each uncountable $\kappa$ with $\cof(\kappa)=\alephes$, $\cof([\kappa]^\alephes)=\kappa^+$.
SSH follows, for example, from the Generalized Continuum Hypothesis.
Shelah's Strong Hypothesis was originally stated differently,
but was shown in \cite[Theorem 6.3]{Sh420} to be equivalent to the variation adopted here.\footnote{In fact,
only the main implication is provided there.
For the other implication: If $\kappa$ is such that $\op{pp}(\kappa)>\kappa^+$, then in particular $\cof[\kappa]^{\cof(\kappa)}>\kappa^+$,
and we may (e.g., by Lemmata 3.4 and 3.8 in \cite{Rinot}) arrange that $\cof(\kappa)=\alephes$.}
The adjective ``Strong'' in SSH means that there is a yet weaker hypothesis, but SSH is in
fact quite weak. In particular, its failure implies the consistency of large cardinals.\footnote{
The failure of SSH at $\kappa$ implies that
in the Dodd--Jensen core model, there is a measurable $\lambda\le\kappa$, moreover $o(\lambda)=\lambda^{++}$.
The exact consistency strength of SSH was established by Gitik in \cite{GitikSCH89, GitikSCH91}.}

Following is the concluding Theorem 6.3 of \cite{Sh420}. The simplicity
of the proof given here is due to the reformulation of SSH.

\bthm[Shelah \cite{Sh420}]\label{SSHthm}
Assume SSH.
For each $\kappa>\alephes$,
the cardinal $\cof([\kappa]^\alephes)$ is $\kappa$ if $\cof(\kappa)>\alephes$, and $\kappa^+$ if $\cof(\kappa)=\alephes$.
\ethm
\bpf
The case $\kappa=\aleph_1$ is Corollary \ref{nojump}. Continue by induction on $\kappa$:
If $\cof(\kappa)=\alephes$, use SSH.
If $\cof(\kappa)>\alephes$, use Lemma \ref{bigcof} and the induction hypothesis to get
$$\cof([\kappa]^\alephes)=\kappa\cdot\sup\set{\cof([\lambda]^\alephes)}{\lambda<\kappa}\le \kappa\cdot
\sup\set{\lambda^+}{\lambda<\kappa}=\kappa.\qedhere$$
\epf

It follows that, assuming SSH, we have that the cardinal $\cof(\kfinw)$ is 
$\fd\cdot\kappa$ if $\cof(\kappa)>\alephes$ and $\fd\cdot\kappa^+$ if $\cof(\kappa)=\alephes$.
Thus, under SSH, the value of $\cof(\kfinw)$ is completely determined.
Moreover, in Theorem \ref{SSHthm}, it suffices to assume that Shelah's Strong Hypothesis holds
for all $\lambda\le\kappa$.

\subsection{Bounds in ZFC}

Even without any hypotheses beyond the ordinary axioms of mathematics,
nontrivial bounds on $\kfinw$ can be established in many cases, using
Shelah's \emph{pcf theory} \cite{ShPCF}.
There are several good introductions to pcf theory. A recent one is
\cite{AbrMagHBK}, whose references include some additional introductions.
The following deep result appears as Theorem 7.2 in \cite{AbrMagHBK}.

\bthm[Shelah]\label{ShCor}
For each $\alpha<\aleph_\alpha$, $\cof([\aleph_\alpha]^{\card{\alpha}})<\aleph_{\card{\alpha}^{+4}}$.
\ethm

In \cite{AbrMagHBK}, Theorem \ref{ShCor}  is stated for limit ordinals $\alpha$, but taking
$\delta=\alpha+\w$, we have that $\delta<\aleph_\alpha<\aleph_\delta$, and applying Shelah's Theorem for
the limit ordinal $\delta$, $\cof([\aleph_\alpha]^{\card{\alpha}})\le\cof([\aleph_{\delta}]^{\card{\alpha}})
=\cof([\aleph_{\delta}]^{\card{\delta}})<\aleph_{\card{\delta}^{+4}}=\aleph_{\card{\alpha}^{+4}}$.


\bdfn
Let \IndexPrint{$\pi$} be the first fixed point of the $\aleph$ function, i.e.,
the first ordinal (necessarily, a cardinal) $\pi$ such that $\pi=\aleph_\pi$.
\edfn

$\pi$ is quite big:
Let $\pi_0=\aleph_0$ and for each $n$, let $\pi_{n+1}=\aleph_{\pi_n}$. Then $\pi=\sup_n\pi_n$.

Shelah's Theorem has the following immediate corollaries.

\bcor\label{lt4}
For each $\alpha<\pi$, $\cof([\aleph_\alpha]^{\alephes})<\aleph_{\card{\alpha}^{+4}}$.
\ecor
\bpf
By induction on $\alpha$.
For $\alpha<\omega$ this follows from Corollary \ref{nojump}.
Assume that the assertion is true for all $\beta<\alpha$, and prove
it for $\alpha$. First, $\cof([\aleph_\alpha]^\alephes)\le \cof([\aleph_{\alpha}]^{\card{\alpha}})\cdot\cof([\card{\alpha}]^\alephes)$.
As $\alpha<\pi$, Corollary \ref{ShCor} is applicable, and thus $\cof([\aleph_{\alpha}]^{\card{\alpha}})<\aleph_{\card{\alpha}^{+4}}$.
Let $\beta$ be such that $\card{\alpha}=\aleph_\beta$. Then $\beta<\pi$,
and thus $\beta<\aleph_\beta=\card{\alpha}$. By the induction hypothesis,
$\cof([\aleph_\beta]^\alephes)<\aleph_{\card{\beta}^{+4}}\le\aleph_{\card{\alpha}^{+3}}$.
\epf

\bcor\label{lt3a}
For each successor cardinal $\kappa<\pi$ and each $\alpha$ with $\kappa\le\alpha<\kappa+\omega$,
we have that
$\cof([\aleph_\alpha]^{\alephes})<\aleph_{\kappa^{+3}}$.
\ecor
\bpf
For each $\beta\in\{\kappa,\kappa+1,\kappa+2,\dots\}$, either $\beta=\kappa$ and
$\cof(\aleph_\beta)=\cof(\kappa)>\alephes$, or $\beta$ is a successor ordinal, and thus
$\cof(\aleph_\beta)=\aleph_\beta>\alephes$. Thus, by Lemma \ref{bigcof},
\begin{eqnarray*}
\cof([\aleph_\alpha]^\alephes) & = & \aleph_\alpha\cdot\sup
\set{\cof([\aleph_\beta]^\alephes)}{\aleph_\beta\le\aleph_{\alpha}, \cof(\aleph_\beta)=\alephes} =\\
& = & \aleph_\alpha\cdot\sup
\set{\cof([\aleph_\beta]^\alephes)}{\beta<\kappa, \cof(\beta)=\alephes} \le\\
& \le & \aleph_\alpha\cdot\sup\set{\cof([\aleph_\beta]^\alephes)}{\beta<\kappa}.
\end{eqnarray*}
By Corollary \ref{lt4}, for each $\beta<\kappa$, $\cof([\aleph_\beta]^\alephes)<\aleph_{\card{\beta}^{+4}}$.

$\aleph_\alpha<\aleph_{\card{\alpha}^+}=\aleph_{\kappa^+}<\aleph_{\kappa^{+3}}$.
Now, for each $\beta<\kappa$, $\cof([\aleph_\beta]^\alephes)<\aleph_{\card{\beta}^{+4}}\le \aleph_{\kappa^{+3}}$.
As $\cof(\aleph_{\kappa^{+3}})=\kappa^{+3}>\kappa$, the supremum is also smaller than
$\aleph_{\kappa^{+3}}$.
\epf

\bcor\label{lt3b}
For each cardinal $\kappa$ with $\alephes<\cof(\kappa)<\kappa<\pi$
and each $\alpha$ with $\kappa\le\alpha<\kappa+\omega$, we have that
$\cof([\aleph_\alpha]^{\alephes})=\aleph_{\alpha}$.
\ecor
\bpf
Replace, in the proof of Corollary \ref{lt3a}, the last paragraph with the following one:
For each $\beta<\kappa$, $\card{\beta}^{+4}<\kappa$, and thus $\aleph_{\card{\beta}^{+4}}<\aleph_{\kappa}\le\aleph_\alpha$.
\epf

\bexm
For each $n\ge 1$:
\be
\itm For each $\alpha<\w_n+\w$, $\cof([\aleph_{\alpha}]^\alephes)<\aleph_{\w_{n+3}}$.
\itm $\cof([\aleph_{\aleph_{\w_n}}]^\alephes)=\aleph_{\aleph_{\w_n}}$.
\ee
\eexm

Combining Theorem \ref{morph} and the estimations provided here for
$\cof([\kappa]^{\alephes})$, we obtain estimations for $\cof(\kfinw)$.

\ed